\documentclass[preprint,1pt]{elsarticle}
\usepackage{pgfplots}
\pgfplotsset{compat=1.17} 

\usepackage{algorithm}
\usepackage{algorithmic}
\usepackage{changepage}
\usepackage[caption=false]{subfig}
\usepackage{multirow}
\usepackage[percent]{overpic}
\usepackage{enumitem}
\newcommand{\customitem}[1]{
    \noindent\hangindent=1.5em\hangafter=1
    \makebox[1em][l]{\textbullet}\ignorespaces #1\par
}

\usepackage{bm}
\usepackage{color}

\newcommand{\blue}[1]{{\color{blue} #1}}

\usepackage{soul}
\usepackage[normalem]{ulem}

\patchcmd\newpage{\vfil}{}{}{}
\usepackage{amssymb}
\usepackage{amsmath}
\allowdisplaybreaks
\usepackage{hyperref}
\usepackage[capitalise]{cleveref}
\journal{review}

\begin{document}
\begin{frontmatter}
\title{Energy-conserving Kansa methods for Hamiltonian wave equations} 
\author[l1]{Xiaobin Li}
\ead{18296320025@163.com}
\affiliation[l1]{organization={School of Mathematics and Computer Sciences},
             addressline={Nanchang University},
             city={Nanchang},
             country={China}}
\author[l1,l2]{Meng Chen}
\ead{chenmeng\_math@ncu.edu.cn}

\affiliation[l2]{organization={Institute of Mathematics and Interdisciplinary Sciences},
             addressline={Nanchang University},
             city={Nanchang},
             country={China}}
\author[l3]{Zhengjie Sun}
\ead{zhengjiesun@njust.edu.cn}
\affiliation[l3]{organization={School of Mathematics and Statistics},
             addressline={Nanjing University of Science and Technology},
             city={Nanjing},
             country={China}}
\author[l4]{Leevan Ling}
\ead{lling@hkbu.edu.hk}
\affiliation[l4]{organization={Department of Mathematics},
             addressline={Hong Kong Baptist University},
             city={Kowloon Tong},
             country={Hong Kong}}
\author[l5]{Siqing Li}
\affiliation[l5]{organization={College of Mathematics},
             addressline={Taiyuan University of Technology},
             city={Taiyuan},
             country={China}}

\ead{lisiqing@tyut.edu.cn}

\begin{abstract}
We introduce a fast, constrained meshfree solver  designed specifically to inherit energy conservation (EC) in second-order time-dependent Hamiltonian wave equations. For discretization, we adopt the Kansa method, also known as the kernel-based collocation method, combined with time-stepping. This approach ensures that the critical structural feature of energy conservation is maintained over time by embedding a quadratic constraint into the definition of the numerical solution.
To address the computational challenges posed by the nonlinearity in the Hamiltonian wave equations and the EC constraint, we propose a fast iterative solver based on the Newton method with successive linearization. This novel solver significantly accelerates the computation, making the method highly effective for practical applications. Numerical comparisons with the traditional secant methods highlight the competitive performance of our scheme. These results demonstrate that our method not only conserves the energy but also offers a promising new direction for solving Hamiltonian wave equations more efficiently.
While we focus on the Kansa method and corresponding convergence theories in this study, the proposed solver is based solely on linear algebra techniques and has the potential to be applied to EC constrained optimization problems arising from other PDE discretization methods.
\end{abstract}

\begin{keyword}

 Energy conservation \sep Hamiltonian wave equations \sep kernel-based collocation methods

 \MSC 65M60 \sep 65M12 \sep 65D30 \sep 41A30
\end{keyword}

\end{frontmatter}

\section{Introduction}
Hamiltonian wave equations are fundamental in describing various natural wave phenomena, including sound waves, electromagnetic waves, and elastic waves \cite{chew2009integral, someda2006electromagnetic}. In this paper, we consider the following second-order Hamiltonian wave equations in a bounded domain $\Omega \subset \mathbb{R}^d$ subject to a boundary condition on $\Gamma=\partial\Omega$:
\begin{subequations}\label{eq:semilinear Hamiltonian wave system}
\begin{align}
	&\ddot u(x,t)-\Delta u(x,t)+F^{'}(u(x,t))=0,  &(x,t)\in \Omega\times\left(0,T\right],
	\label{eq:semilinear Hamiltonian wave equation}
    \\
    & \mathcal{B}u(x,t)= g(x,t),   &(x,t)\in\Gamma\times\left(0,T\right],
    \label{Boundary condition}
    \\
    &u(x,0)=\psi_0(x), \quad \dot u(x,0)=\psi_1(x), \quad & x\in \Omega.
    \label{initial conditions}
\end{align}
\end{subequations}
The term $F^{'}(u)$ in equation \eqref{eq:semilinear Hamiltonian wave equation} represents a conventional nonlinearity, where $F^{'}:\mathbb{R} \rightarrow \mathbb{R}$ is a smooth real-valued function. This type of nonlinearity is frequently encountered in the literature when dealing with Hamiltonian wave equations, as discussed in \cite{reich2000multi}. The initial conditions \eqref{initial conditions} are given by smooth functions $\psi_0$ and $\psi_1$.
We focus on two types of boundary conditions: the Dirichlet boundary condition $\mathcal{B} = I$, and the Neumann boundary condition $\mathcal{B} = \partial_{\vec{n}}$, with $\vec{n}$ being the outward-pointing normal to the boundary $\Gamma$, such that
\[
\int_{\Gamma} \dfrac{\partial u(x,t)}{\partial \vec{n}} \dot{u}(x,t) \, d\Gamma = 0, \text{\qquad for all } t \in [0, T].
\]
Then, the system in \cref{eq:semilinear Hamiltonian wave system} is energy-conservative with  the energy functional defined by
\begin{equation}\label{energy function}
E[u, \dot{u}](t) := \dfrac{1}{2} \int_{\Omega} \left[ (\dot{u}(x, t))^2 + \vert \nabla u(x, t) \vert^{2} + 2F(u(x, t)) \right] dx,
\end{equation}
which is invariant with respect to time. 
It is crucial to develop numerical methods that are structure-preserving, e.g. EC, in the corresponding semi-discrete forms for the original continuous dynamical systems \cite{christiansen2011topics}.
Various structure-preserving methods for partial differential equations (PDEs) and their applications in nonlinear stability have been discussed in the literature \cite{dahlby2011general,gong2019energy}.
Brugnano et al. introduced Hamiltonian boundary value methods (HBVMs) \cite{brugnano2015reprint}, 
a class of Runge-Kutta methods that not only preserve polynomial Hamiltonians exactly, 
but also maintain generic Hamiltonians up to machine accuracy. 
These methods and their spectral version have been applied to Hamiltonian PDEs in \cite{brugnano2015energy,brugnano2019spectrally,brugnano2016line}.
Specifically, finite difference methods preserving local energy conservation, 
particularly for the 1D sine-Gordon equation through non-polynomial conservation laws, 
have been established in \cite{frasca2025finite}.
Additionally, the average vector field (AVF) method, a systematic approach for solving Hamiltonian PDEs, has been presented in \cite{celledoni2012preserving}.
Recently, Cheng et al. \cite{cheng2020new,cheng2020global} employed the Lagrange multiplier method to ensure energy stability in gradient flow model algorithms.
Moreover, this method was applied in developing an energy-conserving algorithm for the coupled Klein-Gordon-Schr{\"o}dinger system, detailed in \cite{li2023high}.
However, the majority of these methods primarily relied on finite difference \cite{hou2019conservative,liu2018new}, 
finite element methods \cite{cai2019linearized,diaz2009energy} and spectral methods \cite{zhang2020highly}.
Therefore, it holds significant theoretical and practical importance to study the energy-conserving numerical methods for solving nonlinear PDEs, particularly in the context of meshfree methods, which have not been extensively explored in the literature.

This paper is organized as follows. Section \ref{sec:Least_squares_Kansa_method} presents the least-squares Kansa method for discretizing the PDEs and the energy conservation principle, along with a study of the solvability and convergence properties at each time step.
Section \ref{sec:Fast_iterative_algorithm} introduces a fast iterative algorithm for solving nonlinear least-squares problems with quadratic constraint, which is crucial for efficiently solving the discretized system obtained in Section \ref{sec:Least_squares_Kansa_method}.
Section  \ref{sec:Numerical experiments}  demonstrates the accuracy, stability, and energy conservation properties of the proposed solver through numerical simulations.
Conclusions are provided in Section \ref{sec:conclusion}.

\section{Least-squares Kansa method with energy conservation}\label{sec:Least_squares_Kansa_method}

In this section, we introduce the least-squares Kansa method with energy conservation for discretizing Hamiltonian wave equations. The proposed method combines the flexibility and accuracy of the radial basis function (RBF) approximation with the robustness of the convergent least-squares approach \cite{chen2022oversampling}. The Kansa method, pioneered by Kansa in 1990 \cite{kansa1990multiquadrics01, kansa1990multiquadrics02}, is an unsymmetric meshless strong-form collocation method that employs RBFs to solve PDEs. This groundbreaking work has inspired a wealth of research on RBFs and their applications in numerical PDEs \cite{chen2011method, kansa2013numerical, li2015compactly, pang2015space}. In the context of Hamiltonian wave equations, Wu and Zhang \cite{wu2015meshless} proposed two meshless energy conservative methods based on RBF approximation for linear cases, while Sun and Ling \cite{sun2022kernel} recently introduced a meshless conservative Galerkin method that incorporates an appropriate projection operator to ensure global energy preservation. However, the strong-form collocation counterpart is still missing, which will be addressed in this work.

\subsection{Discretization using least-squares Kansa method}

The idea begins by presenting a time-stepping scheme for the Hamiltonian wave equation in \cref{eq:semilinear Hamiltonian wave system} and then focus on the spatial discretization using the least-squares Kansa method  \cite{cheung2018h}.
Firstly, let \(v=\dot{u}\). This allows the equation \eqref{eq:semilinear Hamiltonian wave equation} to be rewritten as a system of first-order PDEs:
\begin{equation}\label{eq:uv_system}
\begin{cases}
\dot{u}=v,\\
\dot{v}=\Delta u-F^{'}(u).
\end{cases}
\end{equation}
We then proceed by applying conventional time-stepping schemes for temporal discretization and subsequently implement spatial discretization to obtain the discrete Hamiltonian system.
The time interval $[0, T]$ is  partitioned into $K$ uniform subintervals with a step size of $\tau = T/K$. For any $1\leq k \leq K$,  the semi-discrete solutions are denoted as:
\begin{equation}\label{semi-discrete solutions}
\begin{cases}
u^{k} \approx u(x,t^{k}), \qquad\text{for $x\in \Omega$},\\
v^{k} \approx \dot u(x,t^{k}),  \qquad\text{for $x\in \Omega$}.
\end{cases}
\end{equation}
The Crank-Nicolson (CN) method is applied to \cref{eq:uv_system} as an example, viz.
\begin{equation}\label{Crank-Nicolson integrator}
    \begin{cases}
        &\dfrac{u^{k} - u^{k-1}}{\tau} - \dfrac{v^{k} + v^{k-1}}{2} =0,
        \\[5pt]
        &\dfrac{v^{k} - v^{k-1}}{\tau} - \dfrac{\Delta u^{k} + \Delta u^{k-1}}{2} + F^{'}\left(\dfrac{u^{k} + u^{k-1}}{2}\right) = 0.
    \end{cases}
\end{equation}

For more compactness in notation for nonlinear terms, we make a standard assumption that $F'(0) = 0$, which means that the point $u = 0$ is an equilibrium point of the system. Rewriting the two-step scheme defined by \cref{Crank-Nicolson integrator} with leading indices $k - 1$ and $k-2$, and expressing the resulting  four equations in matrix form yields:
\[
\setlength{\arraycolsep}{3.5pt} 
\renewcommand{\arraystretch}{1.25} 
\begin{bmatrix}
\frac{1}{\tau}   & -\frac{1}{2} & -\frac{1}{2} & 0 \\
0 &  \frac{1}{\tau} & -\frac{1}{\tau} & 0 \\
0 &  0 & -\frac{1}{2} & -\frac{1}{2} \\
0 &  0 & \frac{1}{\tau} & -\frac{1}{\tau}
\end{bmatrix}
\begin{bmatrix}
u^{k} \\
v^{k} \\
v^{k-1} \\
v^{k-2}
\end{bmatrix} =
\begin{bmatrix}
0 & \frac{1}{\tau} & 0  \\
\frac{\Delta}{2}  & \frac{\Delta}{2}  & 0  \\
0 &  -\frac{1}{\tau} & \frac{1}{\tau} \\
0 & \frac{\Delta}{2} & \frac{\Delta}{2}
\end{bmatrix}
\begin{bmatrix}
 u^{k} \\
u^{k-1} \\
u^{k-2} \\
\end{bmatrix}
-F'\left(
\begin{bmatrix}
0 & 0 & 0\\
\frac{1}{2} & \frac{1}{2} & 0\\
0 & 0 & 0\\
0 & \frac{1}{2} & \frac{1}{2}
\end{bmatrix}
\begin{bmatrix}
u^{k}   \\
u^{k-1} \\
u^{k-2}
\end{bmatrix}
\right).
\renewcommand{\arraystretch}{1.2} 
\]
The $4\times 4$ matrix on the left has a determinant of $1/\tau^3>0$, so the above system is always solvable. 
The solutions $u^{k}$ are maintained and expressed as the sum of a linear term with a differential operator and a nonlinear term related to $F'$, as follows:
\begin{equation}\label{cn mtx}
\text{CN:}\qquad\qquad
u^{k} =
{\mathcal{L}}_\Delta\left(
\begin{bmatrix}
 u^{k} \\
u^{k-1} \\
u^{k-2}
\end{bmatrix}
\right)   +
{\mathcal{N}}_{F'}\left(
\begin{bmatrix}
u^{k} \\
u^{k-1} \\
u^{k-2}
\end{bmatrix}
\right), \qquad k\geq 2.
\end{equation}
Note that $\mathcal{L}_\Delta$ and $\mathcal{N}_{F'}$ will be used to denote the generic linear differential operator and nonlinear functions, respectively, throughout the paper. 

Different time discretization schemes can be unified under the same framework. Consequently, the Crank-Nicolson Adams-Bashforth (CNAB) method \cite{ruuth1995implicit} is also considered in this paper.
 In this approach, the CN method is used for the $v$ and $\Delta u$ terms, and the  second-order Adams-Bashforth scheme is applied to the $F^{'}(u)$ term, which leads to
\begin{align}\label{eq:CNAB_scheme}
    \begin{cases}
        &\dfrac{u^{k} - u^{k-1}}{\tau} - \dfrac{v^{k} + v^{k-1}}{2} =0,
        \\[5pt]
        &\dfrac{v^{k} - v^{k-1}}{\tau} - \dfrac{\Delta u^{k} + \Delta u^{k-1}}{2} + F^{'}\left(\dfrac{3u^{k-1} - u^{k-2}}{2}\right) = 0.
    \end{cases}
\end{align}
The CNAB method can be expressed as a linearized (with respect to the absence of $u^{k}$ in ${\mathcal{N}}_{F'}$) four-step method in the following form:
\begin{equation}\label{cnab mtx}
\text{CNAB:}\qquad\qquad
u^{k} =
{\mathcal{L}}_\Delta\left(
\begin{bmatrix}
 u^{k} \\
u^{k-1} \\
u^{k-2}
\end{bmatrix}
\right)   +
{\mathcal{N}}_{F'}\left(
\begin{bmatrix}
u^{k-1} \\
u^{k-2} \\
u^{k-3} \\
\end{bmatrix}
\right), \qquad k\geq 3.
\end{equation}
To sum up, Eqs. \eqref{cn mtx} and \eqref{cnab mtx} can be simplified into
\begin{equation}\label{both eq}
\text{CN \& CNAB:}\qquad\qquad
  u^{k} = \mathcal{L}_\Delta\left(  u^{k}\right) + \mathcal{N}_{F'}\left(u^{k}\right),
\end{equation}
for both the governing equations for $u^{k}$ using 
Eqs. \eqref{cn mtx} or \eqref{cnab mtx}, {see \ref{sec:appendix A1}}.

We are now ready to fully discretize the wave equation 
\eqref{eq:semilinear Hamiltonian wave equation}. 
In this work, we consider symmetric positive definite kernels 
\(\Phi:\mathbb{R}^d \times \mathbb{R}^d \rightarrow \mathbb{R}\) defined via the radial basis functions. We further assume that \(\Phi(\|x-z\|)\) is translation-invariant, with its Fourier transform exhibiting algebraic decay,
\begin{equation}\label{kernels}
c_1(1+\|\omega\|_2^2)^{-m} \leq \widehat{\Phi}(\omega) \leq c_2(1+\|\omega\|_2^2)^{-m},
\end{equation}
for some smoothness order \(m > d/2\).

The family of Wendland's compactly supported kernels \cite{wendland1995piecewise} and the Whittle--Mat\'ern--Sobolev kernels \cite{Matérn1986} are typical examples that satisfy \cref{kernels}.
In particular, the latter is defined as
\begin{equation}\label{Whittle-Matern-def}
    \Phi(x,z) = \,\|x-z\|_2^{\,m-d/2}\,\mathcal{K}_{m-d/2}\bigl(\,\|x-z\|_2\bigr), \quad m > d/2,\, 
\end{equation}
where \(\mathcal{K}_{m-d/2}\) denotes the modified Bessel function of the second kind.

When using \cref{Whittle-Matern-def} to solve second-order strong elliptic PDEs via strong-form kernel-based collocation methods, it has been proven in \cite{cheung2018h} and \cite{schaback2016all2} that the optimal solutions converge in both the $\ell^\infty$ and $\ell^2$ norms; see the original papers for the error estimates.

Let the sets of quasi-uniform trial centers be \(Z = \{z_1, \ldots, z_{n_Z}\} \subset \Omega\), and collocation points be \(X = \{x_1, \ldots, x_{n_X}\} \subset \overline{\Omega}\) and \(Y = \{y_1, \ldots, y_{n_Y}\} \subset \Gamma\).
The fill distance \( h \) in \(\Omega\) with respect to \(Z\) is defined as:
\[
h = h_Z:= \sup_{\xi \in \Omega} \min_{z_j \in Z} \|\xi - z_j\|_{\ell^2(\mathbb{R}^d)}.
\]
Additionally, the assumption \(n_X + n_Y > n_Z\) is made to permit oversampling within the least-squares collocation framework \cite{chen2023kernel,li2019discrete}.

Given a set of trials $Z$ and the kernel $\Phi$, the finite-dimensional trial space $\mathcal{U}_Z$ is defined as follows:
\begin{equation}\label{trial space}
\mathcal{U}_{Z}:=\text{span} \left\lbrace \Phi(\cdot, z_{j}) \mid z_{j} \in Z\right\rbrace,
\end{equation}
where we seek the fully discrete solution
\begin{equation}\label{rbf exp}
u_{Z}^{k} = \sum\limits_{j=1}^{n_Z} \alpha_{j}^{k} \Phi(\cdot, z_{j}) = \Phi(\cdot, Z)\boldsymbol{\alpha}^{k},
\text{\qquad $k\geq 0$}.
\end{equation}
As an example for readers unfamiliar with Kansa's method and for the sake of introducing notation, let us demonstrate what happens when we collocate the boundary condition $\mathcal{B} u_{Z}^{k} = g(x)$ at the boundary points in $Y$. The resulting $n_Y \times n_Z$ matrix system can be compactly expressed as
\[
  [\mathcal{B}\Phi](Y, Z) \boldsymbol{\alpha}^{k} = g(Y),
\]
whose $ij$-entries are given by $[\mathcal{B}\Phi](y_i, z_j)$.
It is notationally more convenient to define the discrete $\ell^2(Y)$-norm for a continuous function $\varphi$ at $Y$ as
\[
    \| \varphi \|_{\ell^2(Y)}^2 := \sum_{i=1}^{n_Y} |\varphi(y_i)|^2,
\]
then, $\big\| [\mathcal{B}\Phi](Y, Z) \boldsymbol{\alpha}^{k} \big\|_{\ell^2(\mathbb{R}^{n_Y})}$ can be written as  $\| \mathcal{B} u_{Z}^{k} \|_{\ell^2(Y)}$.

Given that all $\boldsymbol{\alpha}^{j}$ for $j < k$ are known, we can collocate the governing equation for $u^k$ in \cref{both eq} at the interior collocation points in $X$. For some \( h \)-dependent weight \( h^{-\theta} \), where \( \theta \geq 0 \), as described in \cite[Thm. 2.7]{cheung2018h}, the unconstrained least-squares Kansa solution is defined as
\begin{equation}\label{unconstrained_LS_Kansa_solution}
  u_{Z, LS}^{k} = \mathop{\arg\inf}_{w \in \mathcal{U_Z}}
  \Big(
  \big\| w - \mathcal{L}_\Delta\left(w\right)- \mathcal{N}_{F'}\left(w\right) \big\|_{\ell^2(X)}^2
  + h^{-2\theta} \big\| \mathcal{B} w -g(Y)\big\|_{\ell^2(Y)}^2
  \Big).
\end{equation}
We can consider the equivalent optimization problem in terms of the coefficient vector
\begin{equation}\label{LS sol}
     \boldsymbol{\alpha}_{LS}^{k} = \mathop{\arg\min}_{\boldsymbol{\eta} \in \mathbb{R}^{n_Z}}
  \big\|A \boldsymbol{\eta} - \mathbf{b}_k(\boldsymbol{\eta})\big\|_{\ell^2(\mathbb{R}^{n_X + n_Y})}^2, 
\end{equation}
where matrix $A$ is of size $(n_X + n_Y) \times n_Z$ (see \cref{A matrix} in \ref{sec:appendix A4}), containing both interior and boundary collocation matrices, and is independent of $k$. The vector function $\mathbf{b}_k$, c.f. \cref{cn mtx}, depends on the source function $F^{'}(u)$, boundary values $\mathcal{B}u$, and  previous $\boldsymbol{\alpha}^{j}$ for $j < k$ (see \cref{b vector} in \ref{sec:appendix A4}).

Note that this unconstrained solution $u_{Z, LS}^{k}$ is theoretically proven to be convergent \cite{cheung2018h}. Next, we impose an additional constraint to restrict the trial space to a subset that fulfills the EC constraint.

\subsection{Energy conservation principle in the discretized system}
The energy \cref{energy function} associated with semi-discrete solutions \cref{semi-discrete solutions} at time $t^k$ is given by
\begin{equation}\label{semi energy}
    {E}[u^{k}, v^{k}](t^{k}) := \int_{\Omega} \left[ \dfrac{1}{2}(v^{k})^2 + \dfrac{1}{2} \vert \nabla u^{k} \vert^{2} + F(u^{k}) \right] dx.
\end{equation}
For either \cref{cn mtx} or \eqref{cnab mtx}, similar to the linear algebra manipulation in the previous section, we can again express $[v^{k}, v^{k-1}]^T$ in terms of $u^{j}$ for $j \leq k$ using $\Delta$ and $F'$. In other words, as before, if we omit all previously known solutions from the notation, we can express
\begin{equation}\label{vk(uk)}
   v^k = v^k(u^k),
\end{equation}
using the generic operators \(\mathcal{L}_\Delta\) and \(\mathcal{N}_{F'}\) defined in \cref{both eq} (see \ref{sec:appendix A2} for details).
This allows us to fully discretize the energy with respect to $\boldsymbol{\alpha}^k$ in \cref{rbf exp}.

Inserting \cref{vk(uk)} into \cref{semi energy} and using some quadrature formula, we can express the fully discretized energy in terms of RBF coefficients as 
\begin{equation}\label{approx E}
{E}[u^k, v^k](t^k) \approx  E_w( \boldsymbol{\alpha}^k ) = \mathcal{Q}(\boldsymbol{\alpha}^k) + \mathcal{N}(\boldsymbol{\alpha}^k),
\end{equation}
where \( \mathcal{Q} \) is a quadratic function (see \cref{E_w:Q}), \( \mathcal{N} \) is a nonlinear function depending on \( F \) (see \cref{E_w:N}), and both \( \mathcal{Q} \) and \( \mathcal{N} \) depend on the quadrature weights \( w \).

Therefore, instead of \cref{LS sol},
we can  define the unknown coefficient vector $\boldsymbol\alpha^{k}$ for the fully discretized solution $u_{Z,ECLS}^k:=\Phi(\cdot,Z)\boldsymbol\alpha^{k}$ to the problem \cref{eq:semilinear Hamiltonian wave system} with EC constraint as follows:
\begin{equation}\label{full_dis_sol}
\boldsymbol\alpha^{k} = \mathop{\arg\min}_{\boldsymbol{\eta} \in \mathbb{R}^{n_Z }}
\Big\{
\big\|A \boldsymbol{\eta} - \mathbf{b}_k(\boldsymbol{\eta})\big\|^2 \; \text{s.t. } \; 
{\big\| B \boldsymbol{\eta} - \mathbf{d} \big\|^2 + \mathcal{N}_{F}(\boldsymbol{\eta})=0}
\Big\},
\end{equation}
where $B$ is the Cholesky-like factor of the quadratic form matrix in $\mathcal{Q}$, and $\mathbf{d}$ depends on previous solutions ( see more details in
Appendix Eqs.\eqref{B matrix} and \eqref{d vector} respectively).

\textbf{Remark:}
We can also apply QR decomposition to matrices \(A\) and \(B\) in \cref{full_dis_sol} to enhance computational efficiency \cite{Hammarling+Lucas-updQRfacleasqprob:08}. This approach enables more stable and efficient matrix operations, which is particularly beneficial for handling large data sets.


\section{Fast iterative algorithm for solving nonlinear least-squares problems with a nonlinear  constraint}\label{sec:Fast_iterative_algorithm}

The optimization problem we derived in \cref{full_dis_sol} is a constrained least-squares problem with two nonlinear terms. First, the vector function $ \mathbf{b}_k(\boldsymbol{\eta}) $ within the objective function introduces nonlinearity. Second, the energy constraint is quadratic in nature with a scalar nonlinear function $ \mathcal{N}_{F}(\boldsymbol{\eta}) $. 
In this section, we propose a fast iterative algorithm to enhance computational efficiency and convergence speed. 
This  algorithm uses the Newton method with successive linearization to cope with the nonlinear components and the quadratic constraint.

Our algorithm is based on the combination of Generalized Singular Value Decomposition (GSVD) and the Lagrange multiplier method with linearization. 
Due to the nonlinearity of both the problem and its constraint, we employ an alternating iteration method.

Firstly, we compute the GSVD of matrices $A, B$ in \cref{full_dis_sol} as
\begin{equation}\label{GSVD}
    \begin{aligned}
    A=UCH^{T},\quad B=VSH^{T},
    \end{aligned}
\end{equation}
where $C = \text{diag}(c_1, \ldots, c_{n_Z})$, $S = \text{diag}(s_1, \ldots, s_{n_Z})$ with nonnegative entries,  $U, V$ are unitary matrices, and $H$ is invertible.


Secondly, linearization is applied to \(\mathbf{b}_k(\boldsymbol{\eta})\) and \(\mathcal{N}_{F}(\boldsymbol{\eta})\), followed by the method of Lagrange multipliers to solve \cref{full_dis_sol}, yielding:
\begin{eqnarray}
  C^T\big(  C H^T \boldsymbol{\eta} - U^T\mathbf{b}_k(\boldsymbol{\eta}^{j-1}) \big) + \lambda S^T\big(S H^T \boldsymbol{\eta} - V^T\mathbf{d} \big) &=& 0,
  \label{LM1}\\
  \big\| S H^T \boldsymbol{\eta} - V^T\mathbf{d} \big\|^2 + \mathcal{N}_{F}(\boldsymbol{\eta}^{j-1})  &=& 0.
  \label{LM2}
\end{eqnarray}
Let $\boldsymbol{\zeta} = H^T \boldsymbol{\eta},$
then  \cref{LM1} becomes  a diagonal system:
\begin{equation}
(C^T C + \lambda S^T S) \boldsymbol{\zeta} = C^T U^T\mathbf{b}_k(\boldsymbol{\eta}^{j-1}) + \lambda S^T V^T \mathbf{d}.
\end{equation}
It allows us to express the  $\boldsymbol\zeta$ explicitly in terms of $\lambda$. 

From $t^{k-1}$ to $t^k$,  we consider an alternating iteration on RBF coefficients and Lagrange multiplier $(\boldsymbol\eta^{j}, \lambda_j)\in \mathbb{R}^{n_Z}\times\mathbb{R}$ at some iterative step $j\geq 1$ with some initial values $(\boldsymbol\eta^{0}, \lambda_0)$ based on the solution at $t^{k-1}$.
We define the iterative solution as 
\begin{equation}\label{eta^j}
  \boldsymbol{\eta}^{j}(\lambda) = H^{-T} (C^TC + \lambda S^TS)^{-1} \big(C^TU^T\mathbf{b}_k(\boldsymbol{\eta}^{j-1})+ \lambda S^TV^T \mathbf{d}\big).
\end{equation}
We can now put \cref{eta^j} into the EC
 \eqref{LM2} to express it as
 a function of $\lambda$:
\begin{equation}\label{Ce}
    \mathcal{C}_E( \lambda ) := \big\| S H^T \boldsymbol{\eta}^j(\lambda) - V^T\mathbf{d} \big\|^2 + \mathcal{N}_{F}(\boldsymbol{\eta}^{j-1}).
\end{equation}
It is analogous to a previously proposed solver for least-squares optimization problems with quadratic inequality constraint \cite{Li+Ling-CollMethCaucProb:19}.
The standard method of  Lagrange multipliers requires $\lambda_j$ to be the root of $\mathcal{C}_E( \lambda_j ) = 0$. 

Thirdly, we propose to apply the Newton method and define
\begin{equation}\label{lam_j}
  \lambda_j := \lambda_{j-1} - \frac{\mathcal{C}_E( \lambda_{j-1} )}{\mathcal{C}_E^{'}( \lambda_{j-1} )}.
\end{equation} 
The iteration alternates between: 

(a) \textbf{RBF coefficient update:} For a fixed \(\lambda_{j-1}\), compute \(\boldsymbol{\eta}^{j}\) via \cref{eta^j}.

(b) \textbf{Multiplier update:} For a fixed \(\boldsymbol{\eta}^{j-1}\) within \(\mathcal{C}_E(\lambda)\), compute the new \(\lambda_j\) via \cref{lam_j}.

We iterate until $|\mathcal{C}_E( \lambda_J )| < \text{tol}$ at some iterative step $J$, for some tolerance $\text{tol}$, and the RBF coefficient vector is $\boldsymbol{\alpha}^{k} = \boldsymbol{\eta}^{J}(\lambda_J)$. We summarize the process of updating from $\boldsymbol{\alpha}^{k-1}$ to $\boldsymbol{\alpha}^{k}$ for $k\geq 1$ in  Algorithm \ref{alg:update_alpha}.

\begin{algorithm}[H]
\caption{The Newton method with successive linearization}
\label{alg:update_alpha}
\begin{algorithmic}[0]
\STATE{Step 1} 
\STATE{\qquad\textbf{Initialization:}} 
\STATE{\hspace*{3em}\customitem {Initial solution and Lagrange multiplier: \((\boldsymbol{\eta}^{0}, \lambda_0) = (\boldsymbol{\alpha}^{k-1}, 0)\),}}
\STATE{\hspace*{3em}\customitem{Tolerance for stopping criterion: \text{tol},}}
\STATE{\hspace*{3em}\customitem{Maximum number of iterations: $j_{\text{max}}$.}}
\STATE{Step 2} 
\STATE{\qquad Compute the GSVD  of matrices $A$ and $B$ as in \cref{GSVD}.}
\STATE{\qquad Compute the 
$\boldsymbol{\eta}^{0}(\lambda_0)$ in \cref{eta^j} and
constraint error
$\mathcal{C}_E^{}(\lambda_0)$ in \cref{Ce}.}
\STATE{Step 3} 
\STATE{\qquad $j=0$.}

\qquad\algorithmicwhile { ($j \leq j_{\text{max}}  \text{ and }  |\mathcal{C}_E(\lambda_{j})|> \text{tol}$)}
\STATE{\hspace*{3em} Let $j=j+1$.}
\STATE{ \hspace*{3em} Evaluate  $\lambda_j$ by  \cref{lam_j}.}
\STATE{\hspace*{3em} Compute the 
constraint error
$\mathcal{C}_E^{}(\lambda_j)$ via \cref{Ce} and $\mathcal{C}_E^{'}(\lambda_{j})$.}
\STATE{ \hspace*{3em} Update $\boldsymbol{\eta}^{j}(\lambda_j)$ in \cref{eta^j},
$\boldsymbol b_k(\boldsymbol{\eta}^{j})$, and $\mathcal{N}_F(\boldsymbol{\eta}^{j})$.}

\qquad\algorithmicend
\STATE{Step 4} 
\STATE{\qquad Compute $\boldsymbol\alpha^{k} := \boldsymbol{\eta}^{j}(\lambda_j)$.}
\end{algorithmic}
\end{algorithm}

\textbf{Remark:}
The solvability, stability, and convergence of Kansa methods can be guaranteed by satisfying certain denseness requirements \cite{chen2022oversampling, cheung2018h, ling2006results, ling2008stable, schaback2016all}. Among them, the least-squares Kansa method has demonstrated remarkable success in solving a wide range of PDEs \cite{Li+Ling-CollMethCaucProb:19,chen2023kernel,chen2020extrinsic,li2019discrete}, showcasing its significant potential for further applications.
However, the theoretical analysis of the least-squares Kansa method with an $F$-dependent energy constraint is currently insufficient. 
A major challenge we face is determining a RBF approximation that satisfies the energy constraint, i.e., $E_w$ in \cref{approx E} equals $E_0$.

\section{Numerical experiments}\label{sec:Numerical experiments}
In this section, we demonstrate the accuracy and efficiency of Algorithm \ref{alg:update_alpha} using numerical experiments. 
All methods are defined as follows:
\begin{itemize}
    \item\textbf{EC-LS Kansa-CN}: 
    The energy-conserving least-squares Kansa method with CN scheme for solving \cref{full_dis_sol}.
    \item\textbf{EC-LS Kansa-CNAB}: 
    The energy-conserving least-squares Kansa method with CNAB scheme for solving \cref{full_dis_sol}.
    \item\textbf{MGAVF} \cite{sun2022kernel}: A meshless Galerkin method with the AVF  time integrator.
\end{itemize}
The same definitions apply to \textbf{LS Kansa-CN} and \textbf{LS Kansa-CNAB} without EC.
All least-squares Kansa methods utilize a boundary scaling factor of $\theta = 1.5$ in \cref{unconstrained_LS_Kansa_solution}  
while consistently employing the Whittle-Matérn-Sobolev kernels in \cref{Whittle-Matern-def} with various smoothness orders $m$ and shape parameters $\epsilon$. 
The following four kinds of Hamiltonian equations are considered in numerical tests.
\begin{enumerate}[label=\textbf{PDE 1.}, labelindent=2em, leftmargin=*, align=left]
\item  The 2D linear wave equation is described by 
\begin{subequations}
\begin{align*}
	&\ddot u(\boldsymbol{x},t)-\Delta u(\boldsymbol{x},t) = 0, 
\hphantom{aaaaaaaaaaaa}
    (\boldsymbol{x},t)\in \left( 0,1\right)^2\times\left(0,100\right],
    \\
    & u(\boldsymbol{x},t) = 1,
    \hphantom{aaaaaaaaaaaaaaaaaaaaaaaaa}(\boldsymbol{x},t)\in\Gamma\times\left(0,100\right],
    \\
    &u(\boldsymbol{x},0) = 1, \quad \dot u(\boldsymbol{x},0) =\sqrt{2}\pi \sin(\pi x) \sin(\pi y), 
    \hphantom{aaaaaaaaii}
    x\in \Omega.
\end{align*}
\end{subequations} 
\end{enumerate}
The exact solution is given by
\[u^*(\boldsymbol{x}, t) = \sin(\sqrt{2}\pi t) \sin(\pi x) \sin(\pi y)+1.\]

\begin{enumerate}[label=\textbf{PDE 2.}, labelindent=2em, leftmargin=*, align=left]
    \item The 1D sine-Gordon equation  is expressed as follows: 
\begin{subequations}
\begin{align*}
&\ddot{u}(x,t) - \Delta u(x,t) + \sin(u) = 0,
\hphantom{aa}
(x,t)\in\left(-20,20\right) \times( 0,15],
    \\
    & u_x(\pm 20, t) = 0, 
\hphantom{aaaaaaaaaaaaaaaaaaaaai}(x,t)\in\Gamma\times\left(0,15\right].
\end{align*}
\end{subequations}
\end{enumerate}
The initial conditions can be derived from the exact solution:
\[
u^{*}(x,t)= 4\tan^{-1}
\left(\dfrac{1}{\zeta}\sinh( \dfrac{\zeta t}{\sqrt{1-\zeta^{2}}})/\cosh( \dfrac{x}{\sqrt{1-\zeta^{2}}})\right)
,\zeta\in(0,1).
\]
As $\zeta\rightarrow 1$, the solution gets steeper. 
Here, we set $\zeta=0.9$.

\begin{enumerate}[label=\textbf{PDE 3.}, labelindent=2em, leftmargin=*, align=left]
    \item  The 2D Klein-Gordon equation is specified by:
\begin{subequations}
\begin{align*}
	&\ddot{u}(\boldsymbol{x}, t) - \Delta u(\boldsymbol{x}, t) + u^3 = 0, 
\hphantom{a}\\&\hphantom{aaaaaaaaaaaaaaaaai}
    (\boldsymbol x,t)\in\left\lbrace (x,y)|x^2+y^2 < 10^2\right\rbrace   \times( 0,7],
    \\
    &  u(\boldsymbol{x}, t) = 0,
    \hphantom{aaaaaaaaaaaaaaaaaaaaaaaaii}
    (\boldsymbol x,t) \in \Gamma\times( 0,7],
    \\
    &u(\boldsymbol{x},0) = 2 \mathrm{sech}(\mathrm{cosh}{(x^2 + y^2)}),\quad \dot{u}(\boldsymbol{x},0) = 0, \hphantom{aaaaaa}\boldsymbol x\in\Omega.
\end{align*}
\end{subequations}
\end{enumerate}
It describes the propagation process of an isolated wave.

\begin{enumerate}[label=\textbf{PDE 4.}, labelindent=2em, leftmargin=*, align=left]
    \item The 2D Klein-Gordon equation  is presented as:
\begin{subequations}
\begin{align*}
	&\ddot{u}(\boldsymbol{x}, t) - \Delta u(\boldsymbol{x}, t) + u^5 = 0,
\hphantom{aaaaa}
    (\boldsymbol x,t)\in (-10,10)^2\times(0,20],
    \\
    &  \dfrac{\partial u}{\partial \vec{n}}(\boldsymbol{x},t) = 0, \hphantom{aaaaaaaaaaaaaaaaaaaaaaai}\:
(\boldsymbol{x},t)\in\Gamma\times(0,20],
    \\
    &u(\boldsymbol{x},0) = 0.5 \mathrm{sech}(\cosh(x)) +  0.5 \mathrm{sech}(\cosh(y)),\\
    &\dot{u}(\boldsymbol{x},0) = 0,\hphantom{aaaaaaaaaaaaaaaaaaaaaaaaaaaaaaaaaaai}\boldsymbol{x} \in \Omega.
\end{align*}
\end{subequations}
\end{enumerate}

To measure the error at
 time $t^k$ and compute the convergence order, we use some sets   $\mathcal{X}$ of quasi-uniform and sufficiently dense sets  of evaluation points to define the relative $L^2$ errors between numerical solution $u_{Z, ECLS}^k$ and the exact solution $u$ by
\[
\|u(\cdot, t^k) -  u_{Z, ECLS}^{k}\|_{L^2(\Omega)} \approx \left(\frac{\sum_{x \in \mathcal{X}} \omega_x |u(x, t^k) -  u_{Z, ECLS}^{k}(x)|^2}{\sum_{x \in \mathcal{X}} \omega_x |u^k(x)|^2}\right)^{1/2}.
\]
 Additionally, the relative energy error is computed as follows:
\[
\mathcal{E}_{\text{energy}} = \big|E_w( \boldsymbol{\alpha}^k ) - E_0\big|/\big|E_0\big|,
\]
where energy is discretized at the point set \( P \) of size $n_P$, with quadrature weights computed via the trapezoidal rule.

\textbf{Remark:} 
Since $u^k_Z\big|_P = \Phi(P,Z)\,\boldsymbol{\alpha}^k$, according to \cref{rbf exp},
the energy computed from the original variables \((u^k_Z,v^k_Z)\) coincides with that computed from the RBF coefficients in \cref{approx E}.
Similarly, the relative energy errors computed from either formulation are identical.
As demonstrated in Table \ref{tab:example1_pde1_tolerance_comparison}, the numerical differences between the two approaches are negligible.
\newline

\textbf{Example 1 (energy error tolerance and number of iterations).}
In this experiment, we investigate the impact of varying energy error tolerances on the performance of the proposed EC-LS Kansa-CN and EC-LS Kansa-CNAB methods for solving PDE 1 and 2. 
Regarding energy constraint, we use the precise energies as follows: $E_0 = \pi^2/4$ for PDE 1 and $E_0 = 16/\sqrt{1 - \zeta^2}$ for PDE 2.
Details on the parameters used in temporal and spatial discretization are provided in Table \ref{tab:simulation_settings12}.

 \begin{table}[tbhp]
    \centering
    \caption{\textbf{(Example 1 for PDEs 1-2)} Temporal and spatial discretization parameter settings.}
\label{tab:simulation_settings12}
    \renewcommand{\arraystretch}{1.2}  
    \setlength{\tabcolsep}{6pt}  
    \begin{tabular}{c  c c c c  c c c}
        \hline
  & T   &  $\tau$ &$n_Z$  & $n_X+n_Y$& $n_P$    &   m &$\epsilon$   \\
   \hline
   PDE 1 & 100  & 0.01 &$11^2$  & $30^2$&$101^2$   &   4 &2  \\
    PDE 2 & 15  & 0.01 &$100$  & $400$ &1601  &   5 & 2  \\
        \hline
    \end{tabular}
    \setlength{\tabcolsep}{6pt}  
    \renewcommand{\arraystretch}{1}
\end{table}

\begin{figure}[!htb]
	\centering	
	{
		\label{fig:example1_pde1_L2}
 	\begin{overpic}[width=0.39\textwidth]{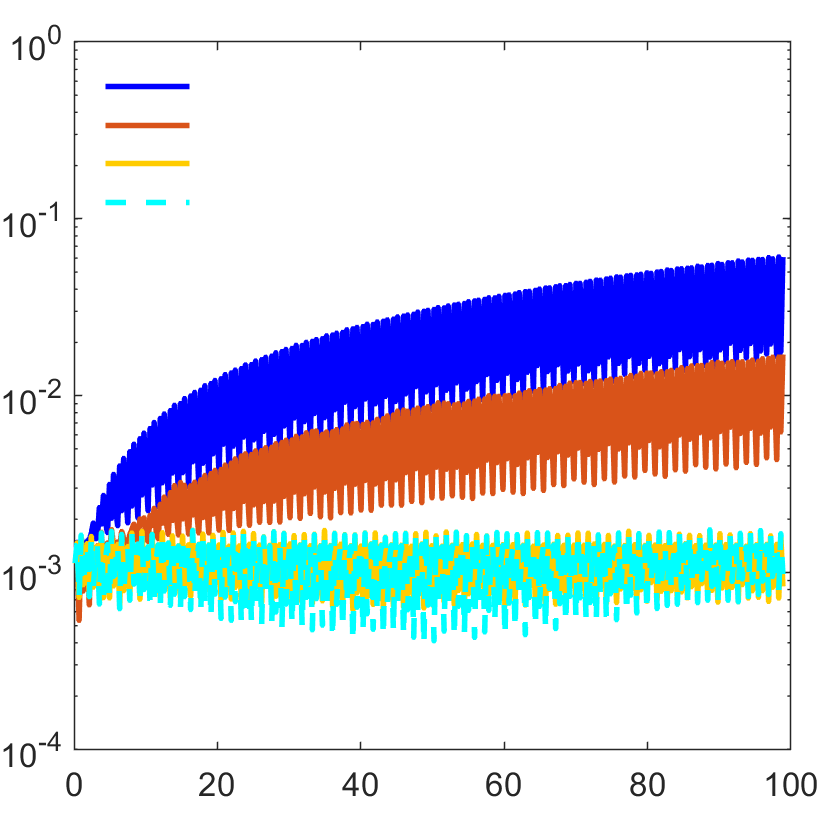}
            \put(-7,42){\scriptsize \rotatebox{90}{$L^{2}\text{ error}$}}
            \put(50,1){\scriptsize \rotatebox{0}{Time}}
            \put(24,88){\scriptsize{$10^{-2}$}}
            \put(24,83.5){\scriptsize{$10^{-4}$}}
            \put(24,79){\scriptsize{$10^{-6}$}}
            \put(24,74.5){\scriptsize{$10^{-8}$}}
 	\end{overpic}
	}\hspace{0pt}
	{
		\label{fig:example1_pde1_energy}
        \begin{overpic}[width=0.39\textwidth]{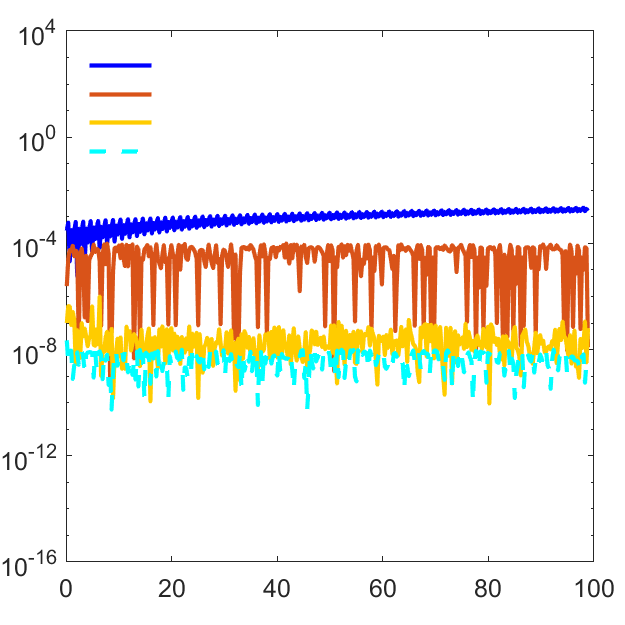}
            \put(-5,45){\scriptsize \rotatebox{90}{$\mathcal{E}_{\text{energy}}$}}
            \put(50,1){\scriptsize \rotatebox{0}{Time}}
           \put(25,88){\scriptsize{$10^{-2}$}}
            \put(25,83.5){\scriptsize{$10^{-4}$}}
            \put(25,79){\scriptsize{$10^{-6}$}}
            \put(25,74.5){\scriptsize{$10^{-8}$}}
 	\end{overpic}
	}
\caption{\textbf{(Example 1 for PDE 1)} The relative $L^2$ errors and relative energy errors vary over time for the EC-LS Kansa-CNAB method at different tolerances, ranging from \(10^{-2}\) to \(10^{-8}\), applied to solving the 2D linear wave equation.}
\label{test3:Figuretol}
\end{figure}
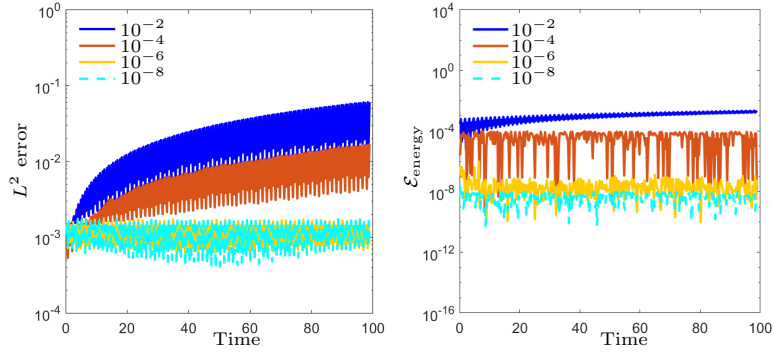 

\begin{table}[!htb]
    \centering
    \caption{\textbf{(Example 1 for PDE 1)} Comparison of relative $L^2$ error, two formulations of relative energy error, iteration counts, and CPU times under different relative energy error tolerances with $\tau = 0.01$ for $T = 100$. Spatial discretization employs a Whittle-Matérn-Sobolev kernel with smoothness order $m=4$ and shape parameter $\epsilon=2$, using $n_Z = 11^2$ and $n_X+n_Y = 30^2$.}
    \label{tab:example1_pde1_tolerance_comparison}
    \renewcommand{\arraystretch}{1.2}  
    \setlength{\tabcolsep}{6pt}  
    \begin{tabular}{c c c c c c}
        \hline
        \multirow{2}{*}{Tolerances} & \multicolumn{4}{c}{EC-LS Kansa-CN(AB)} \\
        & $L^{2}$ errors & $\mathcal{E}_{\text{energy}}(\boldsymbol\alpha^k)$& $\mathcal{E}_{\text{energy}}(u^k_Z,v^k_Z)$ & iterations & CPU(s) \\
        \hline  
        $10^{-2}$ & 1.5630e-02 & 2.0400e-03 & 2.0400e-03 & 10,001 & 1.49 \\
        $10^{-4}$ & 5.4322e-03 & 7.8036e-05 & 7.8036e-05 & 11,178 & 1.54 \\
        $10^{-6}$ & 1.6992e-03 & 9.6382e-09 & 9.2975e-09 & 19,832 & 1.74 \\
        $10^{-8}$ & 1.7185e-03 & 2.7867e-10 & 5.7369e-10 & 52,997 & 1.76 \\
        \hline
    \end{tabular}
    \setlength{\tabcolsep}{6pt}  
    \renewcommand{\arraystretch}{1}  
\end{table}

For PDE 1, the CN and CNAB schemes employ the same formulas due to the linear problem.
As illustrated in Table \ref{tab:example1_pde1_tolerance_comparison}, stricter control of the relative energy error results in a significant reduction in the relative $L^2$ error, although this is accompanied by an increase in the required iteration counts and CPU times. For example, reducing the error tolerance from $10^{-2}$ to $10^{-8}$ leads to an increase in iteration count from 10,001 to 52,997 and CPU time from 1.49 seconds to 1.76 seconds. For further details, refer to Figure \ref{test3:Figuretol}. When the energy error is minimized to a certain threshold, the accuracy of the numerical solution ceases to improve. This finding indicates a trade-off between accuracy and computational cost when adjusting energy error tolerances in numerical simulations. Moreover, Table~\ref{tab:example1_pde1_tolerance_comparison} clearly demonstrates that the energy error computed in terms of the original variables \(u^k_Z\) and \(v^k_Z\) coincides with that computed via the coefficients \(\boldsymbol\alpha^k\).

\begin{table}[!htb]
   \centering
   \caption{\textbf{(Example 1 for PDE 2)} Performance comparison of CN and CNAB schemes under varying relative energy error tolerances, with parameters \(n_Z = 100\), \(n_X+n_Y = 4n_Z\), \(\tau = 0.01\), \(T = 15\), Whittle-Matérn-Sobolev kernel of  \(m = 5\) and  \(\epsilon = 2\).}
    \label{tab:example1_pde2_performance_tolerance_comparison_cn_cnab}
    \renewcommand{\arraystretch}{1.2}  
    \setlength{\tabcolsep}{6pt}  
    \begin{tabular}{c c c c c}
        \hline
        \multirow{2}{*}{Tolerances} & \multicolumn{4}{c}{EC-LS Kansa-CN} \\
        & $L^{2}\text{ errors}$ & $\mathcal{E}_{\text{energy}}$ &  iterations & CPU(s) \\  
        \hline  
        \(10^{-1}\) & 5.7208e-03 & 1.9327e-02 & 1501 & 0.21 \\
        \(10^{-3}\) & 4.3659e-04 & 9.9966e-04 & 1522 & 0.20 \\
        \(10^{-5}\) & 2.9418e-04 & 9.3849e-06 & 4494 & 0.36 \\
        \(10^{-7}\) & 2.9366e-04 & 3.9375e-08 & 4501 & 0.37 \\
        \(10^{-9}\) & 2.9366e-04 & 9.9483e-10 & 4682 & 0.37 \\
        \(10^{-12}\) & 2.9366e-04 & 9.9633e-13 & 7733 & 0.55 \\
        \hline
        \multirow{2}{*}{Tolerances} & \multicolumn{4}{c}{EC-LS Kansa-CNAB} \\
        & $L^{2}\text{ errors}$ & $\mathcal{E}_{\text{energy}}$ & iterations & CPU(s) \\  
        \hline  
        \(10^{-1}\) & 2.9684e-04 & 1.0023e-04 & 1501 & 0.23 \\
        \(10^{-3}\) & 2.9684e-04 & 1.0023e-04 & 1501 & 0.22 \\
        \(10^{-5}\) & 2.9688e-04 & 9.9859e-06 & 1503 & 0.21 \\
        \(10^{-7}\) & 2.9623e-04 & 9.9965e-08 & 1647 & 0.22 \\
        \(10^{-9}\) & 2.9624e-04 & 3.4266e-10 & 3002 & 0.24 \\
        \(10^{-12}\) & 2.9624e-04 & 9.9787e-13 & 3315 & 0.25 \\
    \hline  
       \multirow{2}{*}{Tolerances} & \multicolumn{4}{c}{MGAVF} \\
        & $L^{2}\text{ errors}$ & $\mathcal{E}_{\text{energy}}$ & iterations & CPU(s) \\  
        \hline  
        $10^{-1}$ & 7.3479e-03 & 1.7090e-02 & 1 & 0.17 \\
        $10^{-3}$ & 9.7316e-04 & 9.9714e-04 & 1402 & 0.23 \\
        $10^{-5}$ & 5.3181e-04 & 7.5047e-06 & 1498 & 0.34 \\
        $10^{-7}$ & 5.3296e-04 & 1.0004e-07 & 4923 & 0.35 \\
        $10^{-9}$ & 5.3299e-04 & 1.0676e-09 & 6809 & 0.47 \\  
        $10^{-12}$ & 5.3299e-04 & 1.1911e-10 & 16144 & 1.06 \\  
    \hline
    \end{tabular}
    \setlength{\tabcolsep}{6pt}  
    \renewcommand{\arraystretch}{1}  
\end{table}

Next, we compare the EC-LS Kansa method using the CN and CNAB schemes with the MGAVF method for solving PDE 2, as shown in Table \ref{tab:example1_pde2_performance_tolerance_comparison_cn_cnab}. While the MGAVF method is quite efficient under certain conditions, it often results in higher errors and requires more iterations when the tolerance is tightened. This comparison highlights the EC-LS Kansa method's superior performance, particularly with the CNAB scheme, which not only achieves lower errors but also demands fewer iterations and less CPU time. Importantly, unlike many traditional methods that are tied to a specific time discretization strategy, the EC-LS Kansa method offers the flexibility to use various time-stepping schemes, adapting easily to different computational challenges. Even under stricter error tolerance conditions (less than $10^{-5}$), the performance of the CNAB scheme is comparable to that of the CN scheme, with only minor differences.

  We also compare  EC-LS Kansa-CNAB method with the traditional optimization solver,
  which applies the Lagrange multiplier method and the secant method \blue{\cite{Argyros1988}}. As shown in Table \ref{tbl:example1_timecost}, our method consistently reduces CPU time across all PDEs, with reductions ranging from approximately sixfold to over tenfold. These results underscore the efficiency and potential of our method for solving complex PDEs more rapidly.
  \newline
   \begin{table}[tbhp]
    \centering
    \caption{\textbf{(Example 1 for PDEs 1-4)} Comparison between our proposed method and the traditional secant method to \cref{full_dis_sol}, in solving four PDEs. The final times are \(T = 100\), \(15\), \(7\), and \(20\), with $\tau = 0.01$ for the PDEs 1-3 and  $\tau = 0.05$ for the 4.}
    \renewcommand{\arraystretch}{1.2}  
    \label{tbl:example1_timecost}
    \begin{tabular}{ccccc}
        \hline
        \multirow{2}{*}{Method} & \multicolumn{4}{c}{CPU (s)} \\
        & PDE 1 & PDE 2 & PDE 3 & PDE 4 \\
        \hline
        EC-LS Kansa-CNAB & 1.76 &0.25  &53.07  &66.10  \\
        Secant method    &  13.17    &4.57     &511.12  &451.93  \\
        \hline
    \end{tabular}
    \renewcommand{\arraystretch}{1}  
\end{table}

\textbf{Example 2 (Temporal and spatial  convergence).}
We first test the temporal convergence characteristics of the EC-LS Kansa method by solving PDE 2, using fixed spatial discretization parameters: \(n_Z = 400\), \(n_X + n_Y = 2.5n_Z\), a kernel with \(m = 4\) and \(\epsilon = 4\), along with an energy constraint error tolerance of \(10^{-11}\).

We employ the CN and CNAB schemes to assess the temporal convergence of the EC-LS Kansa method in solving PDEs. We reduce the time step from \(0.04\) to \(0.0025\), and calculate the relative \(L^2\) errors of the numerical solution \(u_{Z,ECLS}^{k}\) at \(T = 1\).

\begin{table}[!htb]
    \centering
    \caption{\textbf{(Example 2 for PDE 2)} Comparison of temporal convergence between different time discretization schemes for the EC-LS Kansa method at $T=1$, with spatial discretization parameters: $n_Z = 400$, $n_X+n_Y = 2.5n_Z$, Whittle-Matérn-Sobolev kernel of  $m = 4$ and $\epsilon = 4$.}
    
    \label{tab:Temporal convergence}
    \renewcommand{\arraystretch}{1.2}  
    \setlength{\tabcolsep}{6pt}  
    \begin{tabular}{c  c c c c  c c c c}
        \hline
        \multirow{2}{*}{$\tau$} & \multicolumn{3}{c}{CN scheme} & \multicolumn{3}{c}{CNAB scheme} \\
         & $L^{2}\text{ errors}$ & rate  & CPU(s) & $L^{2}\text{ errors}$ & rate  & CPU(s) \\
        \hline

0.04    & 3.9246e-04 &        & 0.53   & 4.2593e-04 &        & 0.32 \\
0.02    & 9.8375e-05 & 2.00   & 0.54   & 1.0660e-04 & 2.00   & 0.35 \\
0.01    & 2.4610e-05 & 2.00   & 0.81   & 2.6643e-05 & 2.00   & 0.43 \\
0.005   & 6.1541e-06 & 2.00   & 1.28   & 6.6591e-06 & 2.00   & 0.54 \\
0.0025  & 1.5334e-06 & 2.00   & 1.83   & 1.6590e-06 & 2.00   & 0.81 \\
        \hline
    \end{tabular}
    \setlength{\tabcolsep}{6pt}  
    \renewcommand{\arraystretch}{1}
\end{table}

As shown in Table \ref{tab:Temporal convergence}, with the reduction of the time step \(\tau\), both schemes exhibit the expected second-order convergence in \(L^2\) error, confirming their efficiency in time discretization. In terms of CPU time efficiency, the CNAB scheme consistently outperforms the CN scheme. For example, at $\tau = 0.0025$, the CN scheme requires 1.83 seconds, whereas the CNAB scheme only needs 0.81 seconds. This pattern is observed for all time steps, indicating that the CNAB scheme is not only comparably accurate but also more efficient computationally.
While the CN scheme shows a slight edge in accuracy, the CNAB scheme offers a better balance between accuracy and computational efficiency.

Next, we examine the spatial convergence characteristics under varying settings of kernel smoothness order and oversampling strategies. A sufficiently small time step \(\tau = 5 \times 10^{-4}\) is set, using a kernel with  \(\epsilon = 5\). Center points \(Z\), and collocation points \(X\) and \(Y\) are uniformly selected within the domain \(\Omega = (-20, 20)\).
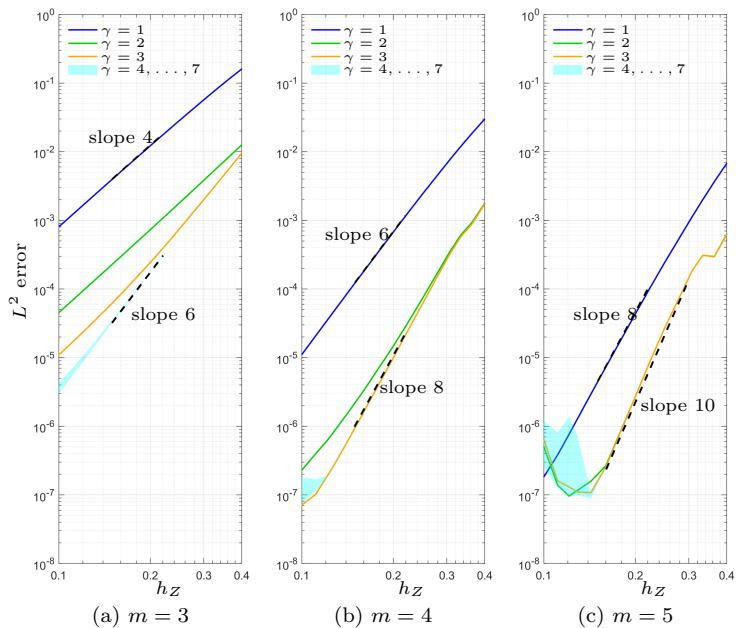
\begin{figure}[!htb]
    \centering
    \subfloat[$m = 3$]
    {
    \begin{overpic}[width=.24\textwidth]{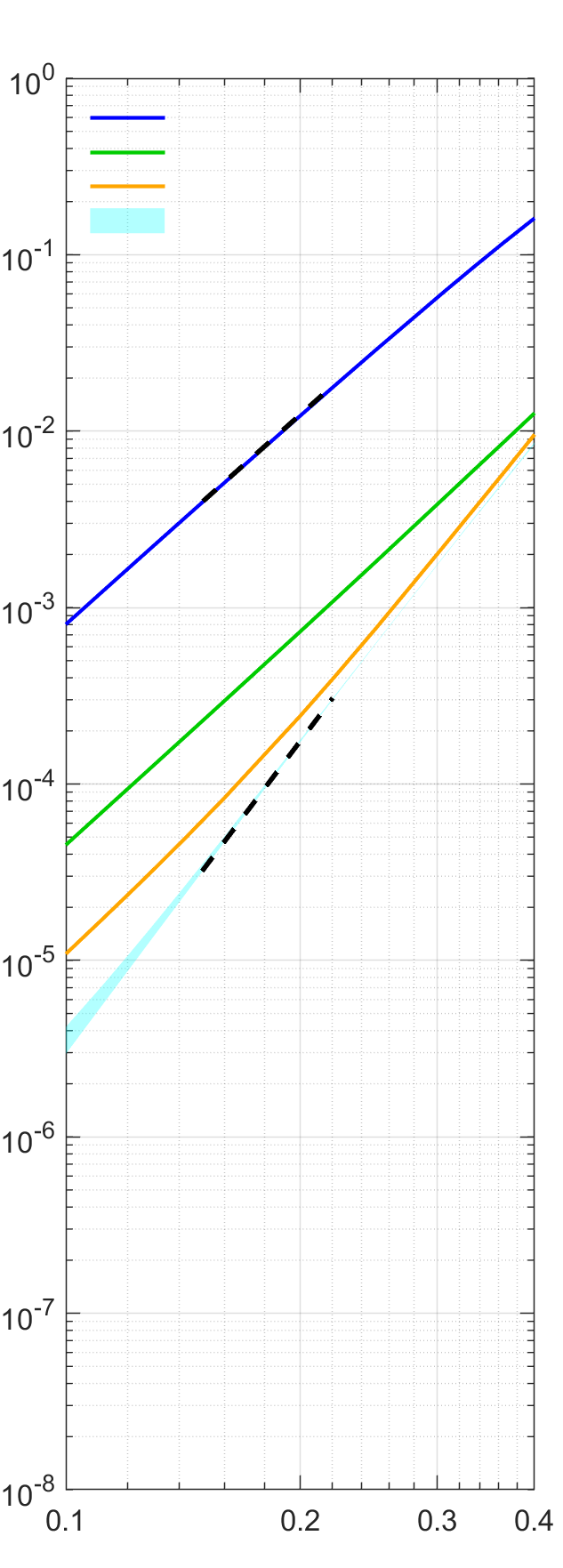}
        \put(-4,45){\scriptsize \rotatebox{90}{$L^2$ error}}
        \put(20,0){\scriptsize \rotatebox{0}{$h_Z$}}
        \put(9,74){\scriptsize \rotatebox{0}{slope 4}} 
        \put(16,45){\scriptsize \rotatebox{0}{slope 6}} 
        \put(11,91.8){\tiny \rotatebox{0}{$\gamma=1$}}
        \put(11,89.7){\tiny \rotatebox{0}{$\gamma=2$}}
        \put(11,87.5){\tiny \rotatebox{0}{$\gamma=3$}}
        \put(11,85.4){\tiny \rotatebox{0}{$\gamma = 4, \ldots, 7$ }}
    \end{overpic}
   }\hspace{1pt}
   \subfloat[$m = 4$]
   {
   \begin{overpic}[width=.24\linewidth]{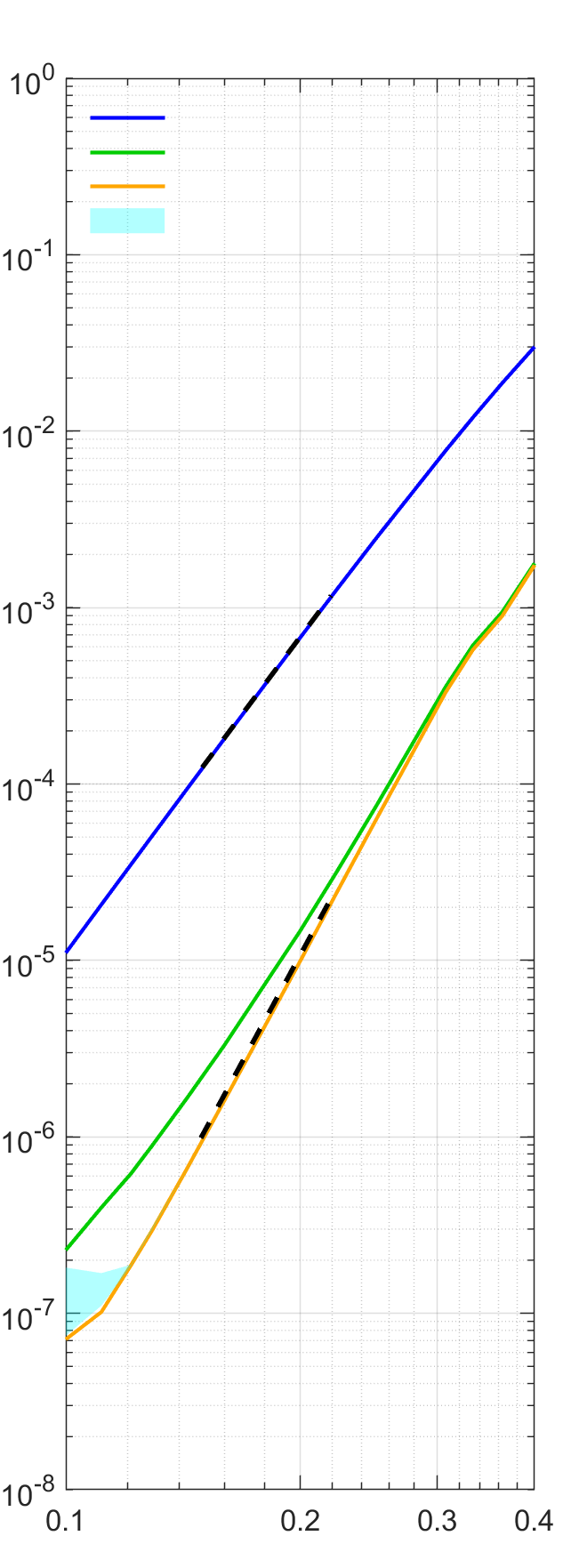}
        \put(20,0){\scriptsize \rotatebox{0}{$h_Z$}}
        \put(8,58){\scriptsize \rotatebox{0}{slope 6}} 
        \put(17,33){\scriptsize \rotatebox{0}{slope 8}} 
        \put(11,91.8){\tiny \rotatebox{0}{$\gamma=1$}}
        \put(11,89.7){\tiny \rotatebox{0}{$\gamma=2$}}
        \put(11,87.5){\tiny \rotatebox{0}{$\gamma=3$}}
        \put(11,85.4){\tiny \rotatebox{0}{$\gamma = 4, \ldots, 7$ }}
    \end{overpic}
    }\hspace{1pt}
    \subfloat[$m = 5$]
   {
   \begin{overpic}[width=.24\linewidth]{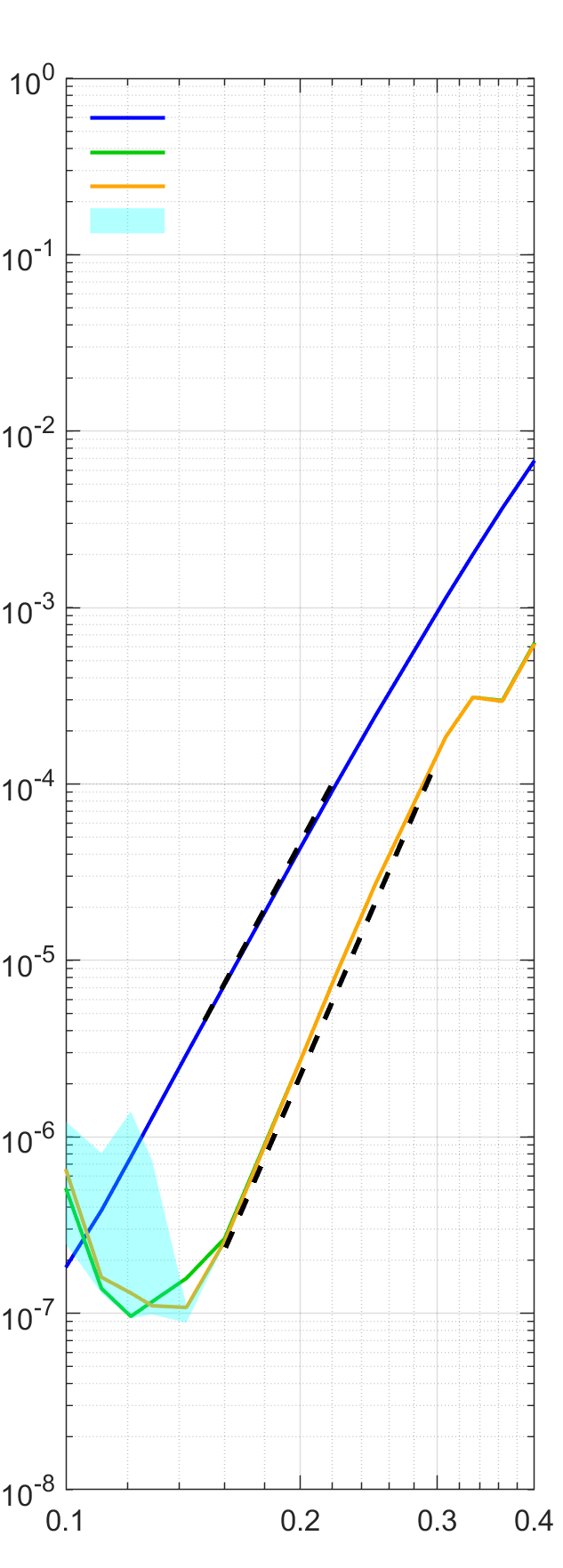}
        \put(20,0){\scriptsize \rotatebox{0}{$h_Z$}}
        \put(9,45){\scriptsize \rotatebox{0}{slope 8}} 
        \put(20,30){\scriptsize \rotatebox{0}{slope 10}} 
        \put(11,91.8){\tiny \rotatebox{0}{$\gamma=1$}}
        \put(11,89.7){\tiny \rotatebox{0}{$\gamma=2$}}
        \put(11,87.5){\tiny \rotatebox{0}{$\gamma=3$}}
        \put(11,85.4){\tiny \rotatebox{0}{$\gamma = 4, \ldots, 7$ }}
   \end{overpic}
   }
    \caption{\textbf{(Example 2 for PDE 2)} Relative \(L^2\) convergence profiles for \textbf{PDE 2} solved using the EC-LS Kansa-CNAB method with Whittle-Matérn-Sobolev kernels of smoothness order \(m = 3, 4, 5\), and \(\epsilon = 5\), for \(\tau = 5\times 10^{-4}\) and \(T = 1\).
    Errors for ratios of oversampling  $ \gamma = (n_X+n_Y)/n_Z, \in \{1, 2, 3\}$ are shown individually, while errors for $ \gamma \in \{4, 5, 6, 7\}$ are collectively represented in a shaded area.}
    \label{fig:spatial  convergence}
\end{figure}

As shown in Figure \ref{fig:spatial  convergence}, the results indicate that spatial convergence rates depend significantly on the smoothness order of the kernel. When \(n_Z = n_X + n_Y\), the convergence order is close to \(2(m-1)\); with higher oversampling rates, this order almost reaches  \(2m\). 
\newline

\textbf{Example 3(Exploring energy conservation on Dirichlet and Neumann boundary conditions).}
In this experiment, we solve PDE 3 and PDE 4 to investigate the energy conservation properties of the EC-LS Kansa-CNAB method under Dirichlet and Neumann boundary conditions respectively. The specific parameter settings are show in Table \ref{tab:simulation_settings}.
 \begin{table}[tbhp]
    \centering
    \caption{\textbf{(Example 3 for PDEs 3-4)} Temporal and spatial discretization parameter settings.}
\label{tab:simulation_settings}
    \renewcommand{\arraystretch}{1.2}  
    \setlength{\tabcolsep}{6pt}  
    \begin{tabular}{c  c c c c  c c c}
        \hline
  & T   &  $\tau$ &$n_Z$  & $n_X+n_Y$& $n_P$    &   $m$ &$\epsilon$   \\
   \hline
   PDE 3 & 7  & 0.01 &2267  & 5809&$101^2$   &   4 &2  \\
    PDE 4 & 20  & 0.05 &$51^2$  & $101^2$ &$101^2$  &   4 &1  \\
        \hline
    \end{tabular}
    \setlength{\tabcolsep}{6pt}  
    \renewcommand{\arraystretch}{1}
\end{table}

\begin{figure}[tbhp]
    \centering	
    \subfloat[Dirichlet boundary]
    {
        \label{fig:example3 Dirichlet boundary}
        \begin{overpic}[width=0.465\textwidth]{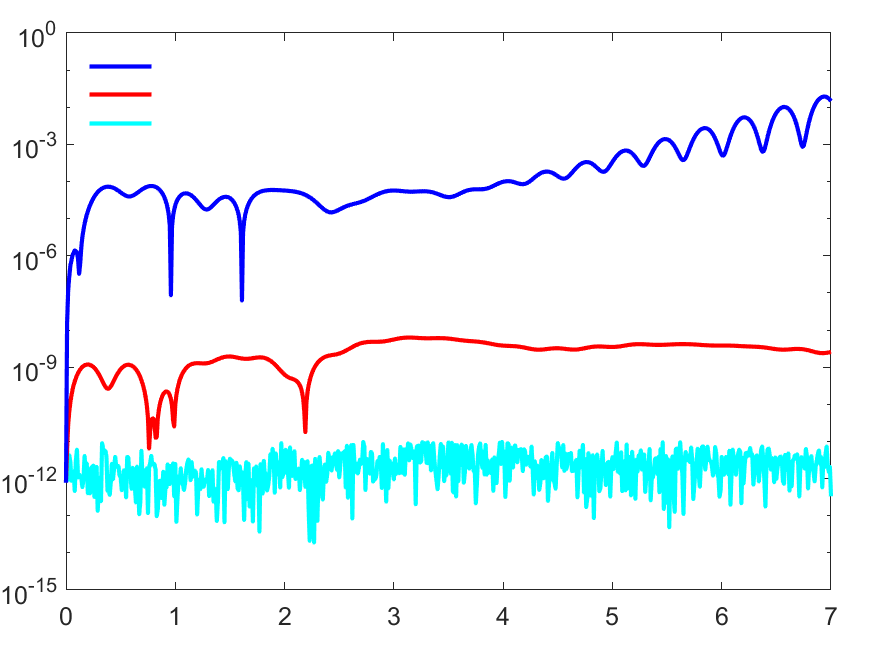}
        \put(-5,32){\scriptsize \rotatebox{90}{$\mathcal{E}_{\text{energy}}$}}
        \put(45,-1){\scriptsize \rotatebox{0}{Time}}
        \put(18,66.5){\tiny \rotatebox{0}{LS Kansa-CNAB}}
        \put(18,63.2){\tiny \rotatebox{0}{MGAVF}}
        \put(18,59.7){\tiny \rotatebox{0}{EC-LS Kansa-CNAB}}
        \end{overpic}
        }\hspace{1pt}
    \subfloat[Neumann boundary]
    {
        \label{fig:example3 Neumann boundary}
        \begin{overpic}[width=0.465\textwidth]{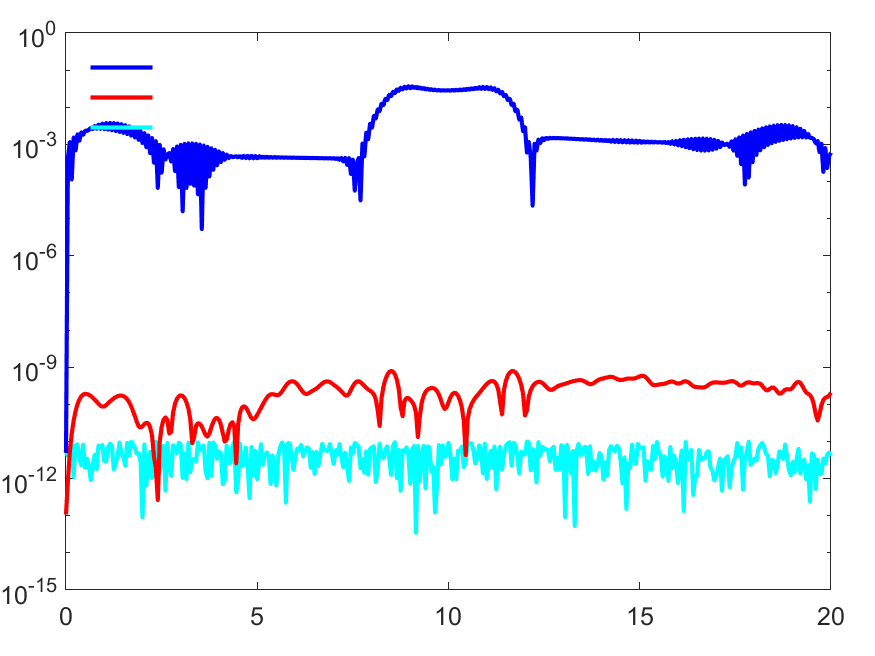}
        \put(-5,32){\scriptsize \rotatebox{90}{$\mathcal{E}_{\text{energy}}$}}
        \put(45,-1){\scriptsize \rotatebox{0}{Time}}
        \put(18,66.5){\tiny \rotatebox{0}{LS Kansa-CNAB}}
        \put(18,63.2){\tiny \rotatebox{0}{MGAVF}}
        \put(18,59.7){\tiny \rotatebox{0}{EC-LS Kansa-CNAB}}
        \end{overpic}
    }
    \caption{\textbf{(Example 3 for PDEs 3-4)} Comparison of energy conservation under different boundary conditions for the 2D Klein-Gordon equation: (a) Relative energy errors for PDE 3 under Dirichlet conditions; (b) Relative energy  errors for PDE 4 under Neumann conditions.}
    \label{fig:example3 different boundary}
\end{figure}
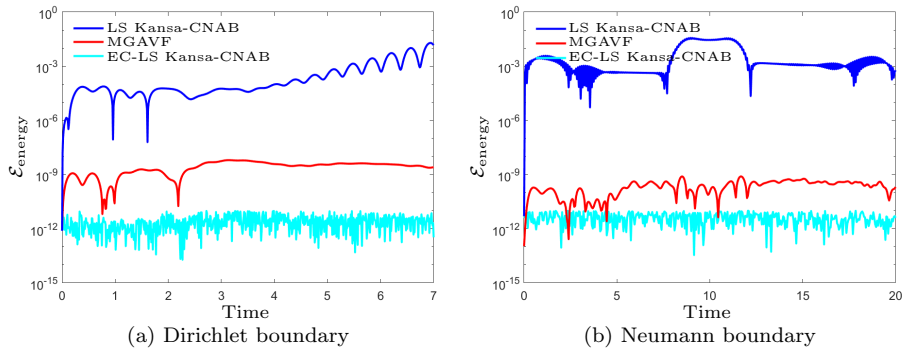

As depicted in Figure \ref{fig:example3 different boundary}, the least-squares Kansa method without EC \cite{chen2023kernel} exhibits poor performance in terms of EC. In contrast, our newly developed EC-LS Kansa-CNAB method achieves EC. Besides, compared to the MGAVF method, our EC-LS Kansa-CNAB method exhibits superior energy conservation.
\newline

\textbf{Example 4 (Numerical accuracy and energy errors under different point distributions).}
In this experiment, we focus on solving PDE 2 with $\tau = 0.01$ for $T = 15$.
The spatial discretization parameters are fixed as follows: \(n_Z = 200\), \(n_X + n_Y = 2n_Z\), kernel of  \(m = 4\) and  \(\epsilon = 2.5\).

\begin{figure}[tbhp]
    \centering
    {
        \label{fig:example4 uniform point}
        \begin{overpic}[width=0.45\textwidth]{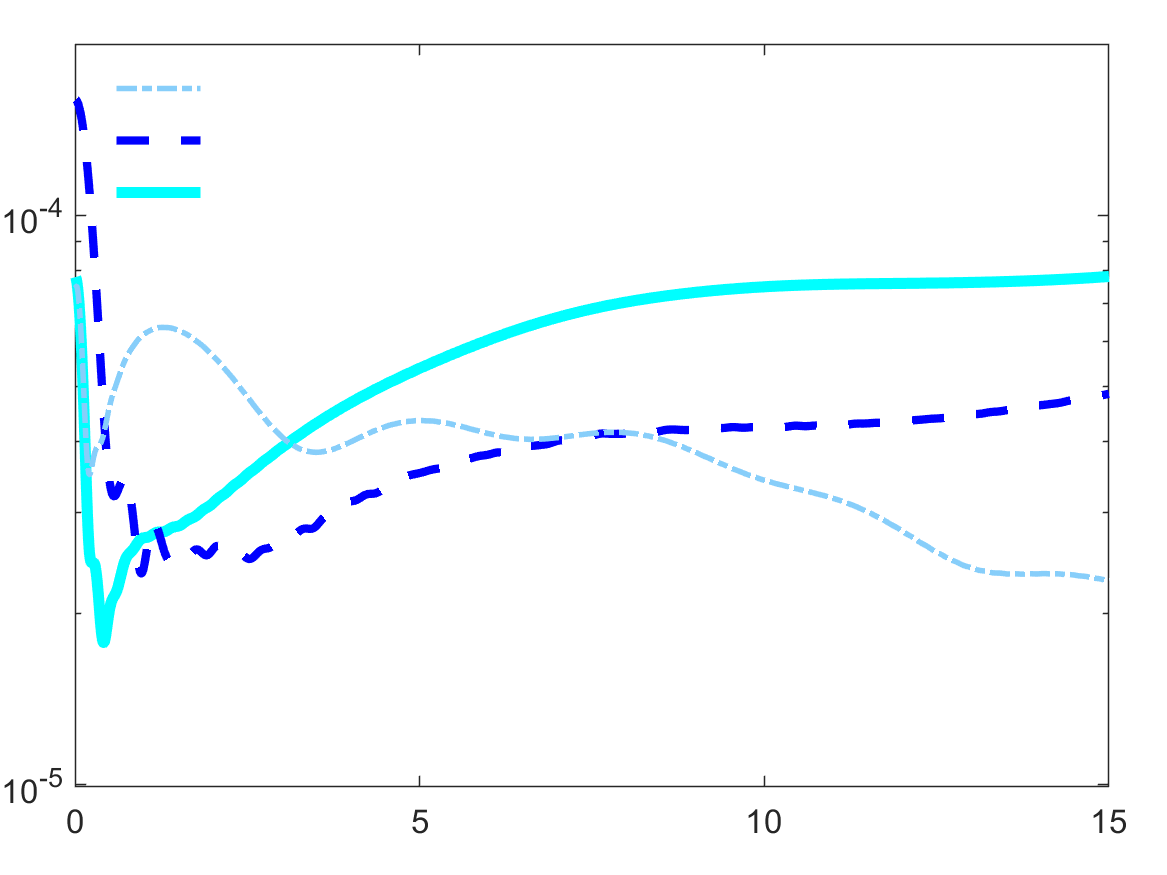}
            \put(-5,30){\scriptsize \rotatebox{90}{$L^2$ error}}
            \put(40,1){\scriptsize \rotatebox{0}{Time}}
            \put(20,66.3){\scriptsize \rotatebox{0}{LS Kansa-CNAB}}
            \put(20,61.9){\scriptsize \rotatebox{0}{MGAVF}}
            \put(20,57.5){\scriptsize \rotatebox{0}{EC-LS Kansa-CNAB}}
        \end{overpic}
    }\hspace{1pt}
    {
        \label{fig:example4 Halton point}
        \begin{overpic}[width=.45\linewidth]{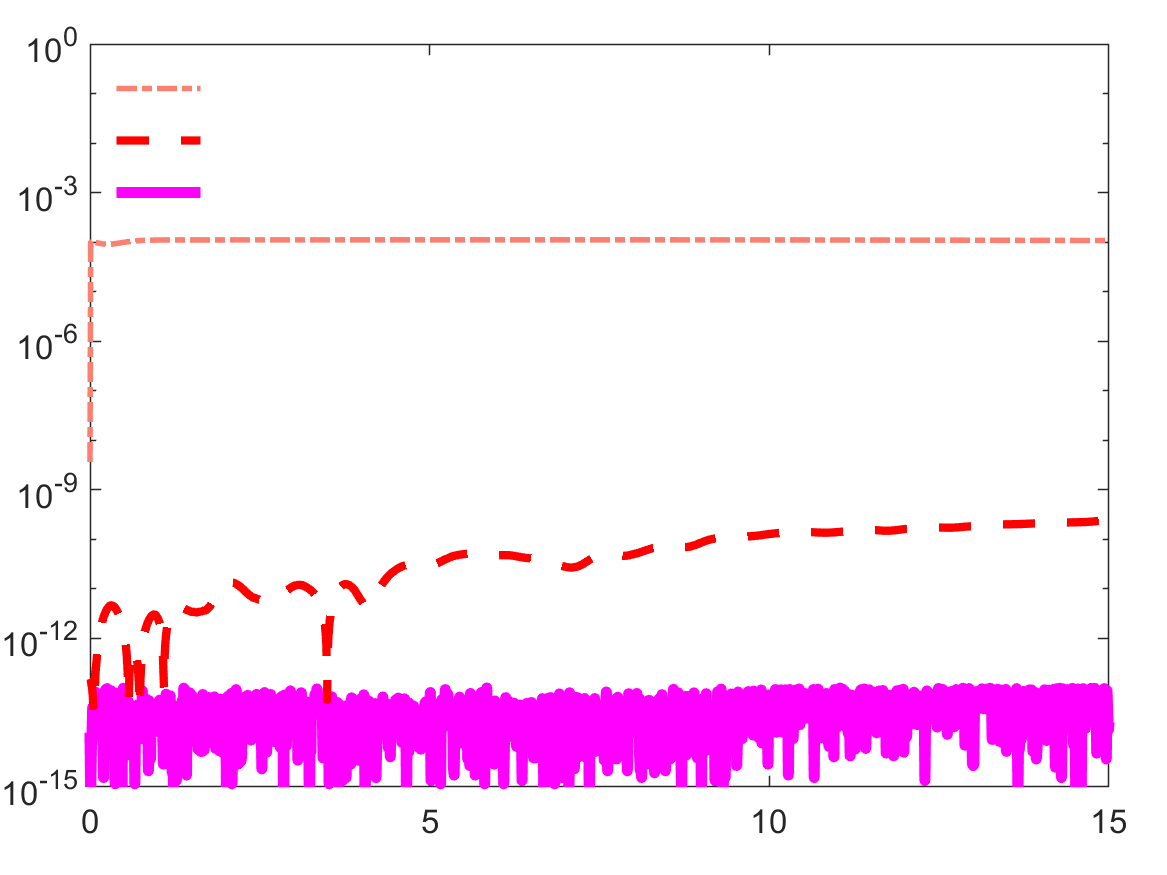}
            \put(-5,32){\scriptsize \rotatebox{90}{$\mathcal{E}_{\text{energy}}$}}
            \put(40,1){\scriptsize \rotatebox{0}{Time}}
            \put(20,66.3){\scriptsize \rotatebox{0}{LS Kansa-CNAB}}
            \put(20,61.9){\scriptsize \rotatebox{0}{MGAVF}}
            \put(20,57.5){\scriptsize \rotatebox{0}{EC-LS Kansa-CNAB}}
        \end{overpic}
    }
    \\
   (a) On uniform points
    \\
    {
        \label{fig:example4 uniform point}
        \begin{overpic}[width=0.45\textwidth]{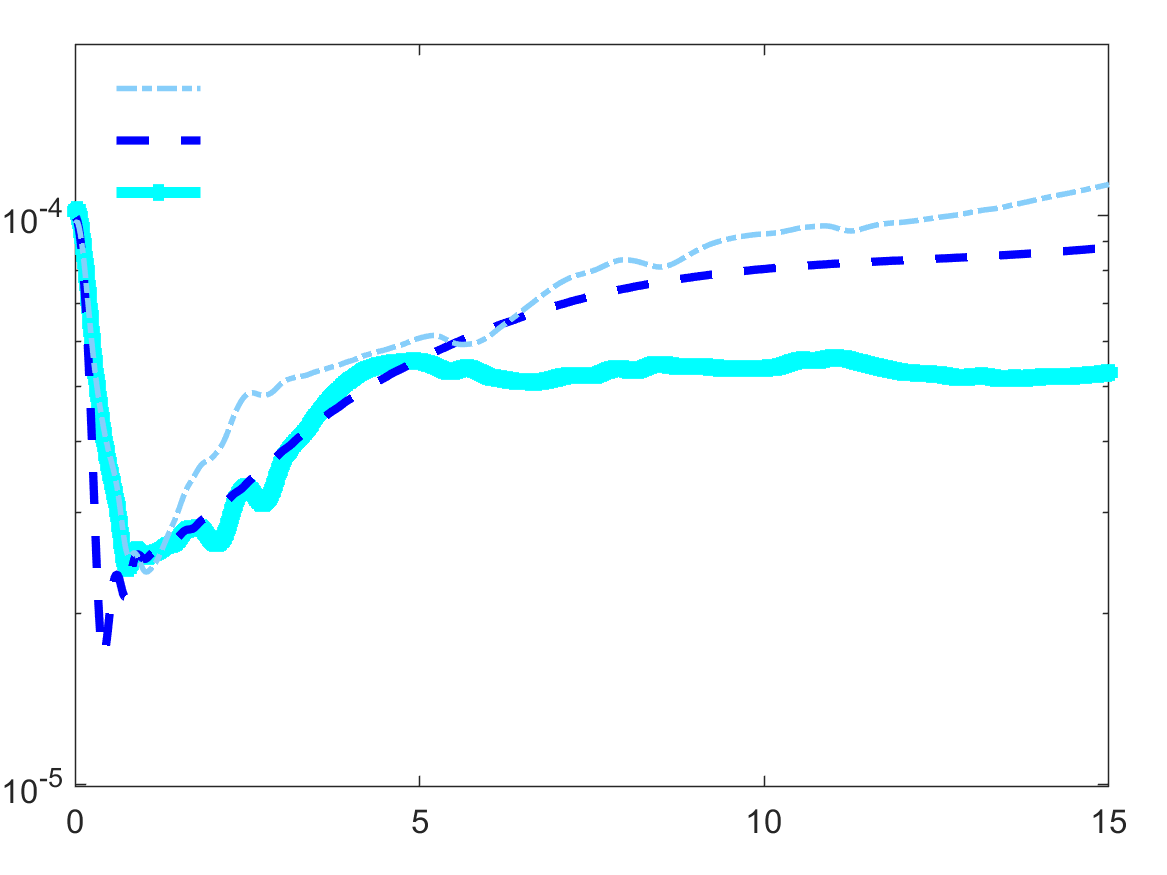}
            \put(-5,30){\scriptsize \rotatebox{90}{$L^2$ error}}
            \put(40,1){\scriptsize \rotatebox{0}{Time}}
            \put(20,66.3){\scriptsize \rotatebox{0}{LS Kansa-CNAB}}
            \put(20,61.9){\scriptsize \rotatebox{0}{MGAVF}}
            \put(20,57.5){\scriptsize \rotatebox{0}{EC-LS Kansa-CNAB}}
        \end{overpic}
    }\hspace{1pt}
    {
        \label{fig:example4 Halton point}
        \begin{overpic}[width=.45\linewidth]{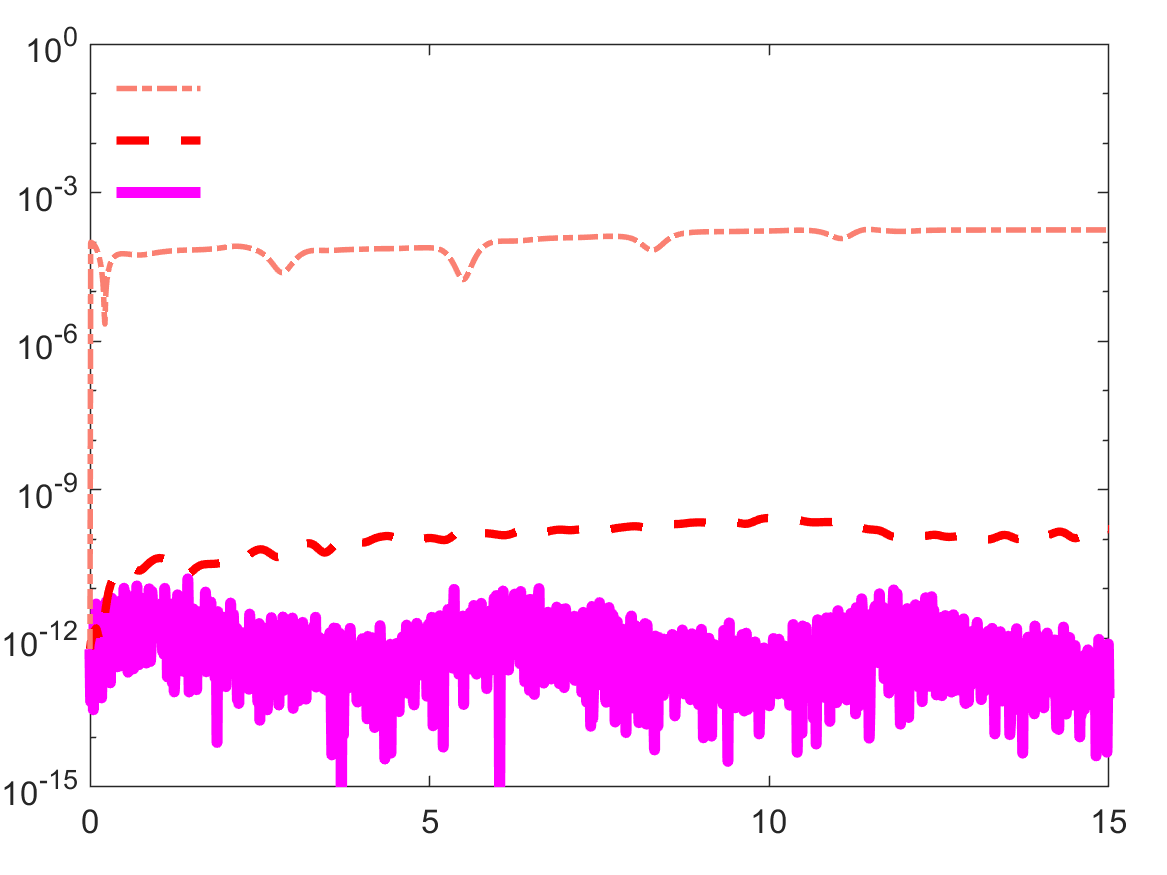}
            \put(-5,32){\scriptsize \rotatebox{90}{$\mathcal{E}_{\text{energy}}$}}
            \put(40,1){\scriptsize \rotatebox{0}{Time}}
            \put(20,66.3){\scriptsize \rotatebox{0}{LS Kansa-CNAB}}
            \put(20,61.9){\scriptsize \rotatebox{0}{MGAVF}}
            \put(20,57.5){\scriptsize \rotatebox{0}{EC-LS Kansa-CNAB}}
        \end{overpic}
    }
    \\
       (b) On Halton points
    \\
    \caption{\textbf{(Example 4 for PDE 2)} Long-time numerical accuracy and energy conservation performance of different methods (LS Kansa-CNAB, MGAVF, and EC-LS Kansa-CNAB) at uniform and Halton point distributions. On (a) uniform points and (b) Halton points, the relative $L^2$ error (cyan, dark blue, and light blue corresponding to EC-LS Kansa-CNAB, MGAVF, and LS Kansa-CNAB) and the associated relative energy error (magenta, red, and light coral corresponding to EC-LS Kansa-CNAB, MGAVF, and LS Kansa-CNAB) are plotted respectively.
} 
    \label{fig:example4}
\end{figure}
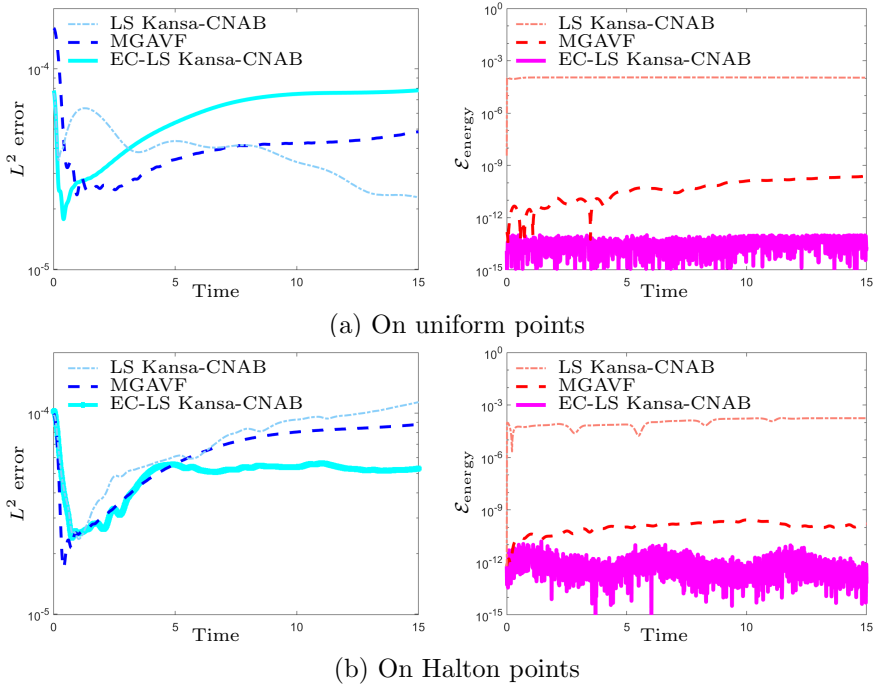

As shown in Figure \ref{fig:example4}, under various point distributions, the EC-LS Kansa-CNAB method demonstrates high numerical accuracy and excellent energy conservation performance in solving PDE 2. In comparison, the performances of EC-LS Kansa-CNAB and MGAVF are less sensitive to point distribution. Particularly, EC-LS Kansa-CNAB surpasses other methods in energy conservation accuracy. Even under Halton point distribution, its accuracy is only slightly affected, maintaining at the $10^{-11}$ magnitude level. Compared to the other two methods, the LS Kansa method is the most sensitive to point distribution, showing the poorest performance in energy conservation.
\newline

\textbf{Example 5 (Long-time numerical simulations of wave equations using the EC-LS Kansa-CNAB method).}
In this example, the EC-LS Kansa-CNAB method is employed in long-time numerical simulations for  PDEs 1, 3, and 4 with different time steps 0.01, 0.01, and 0.05, respectively. 
The detailed parameter settings are presented in Tables \ref{tab:simulation_settings12} and \ref{tab:simulation_settings} respectively.

\begin{figure}[tbhp]
    \centering
        \subfloat[$10^{\mathrm{th}}$ wave period, crest]{\includegraphics[width=.38\linewidth]
        {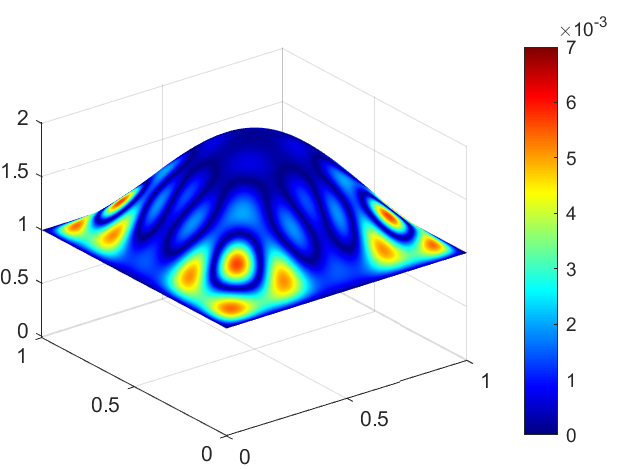}}
        \subfloat[$20^{\mathrm{th}}$ wave period, trough]{\includegraphics[width=.38\linewidth]              
            {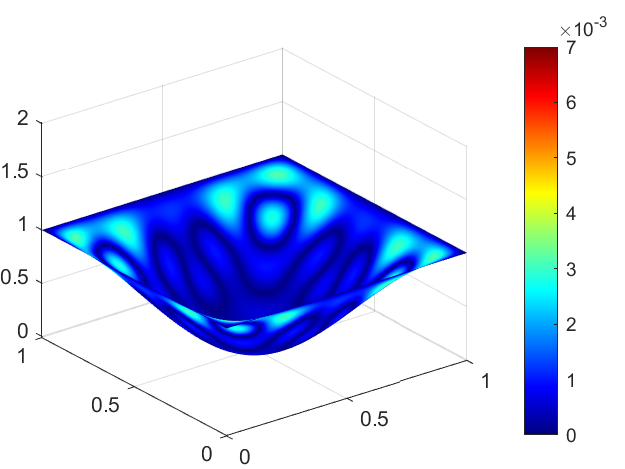}}\\
        \subfloat[$30^{\mathrm{th}}$ wave period, crest]{\includegraphics[width=.38\linewidth]       
             {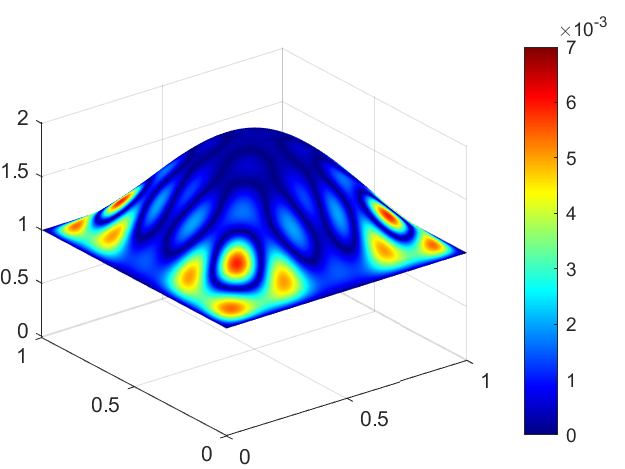}}\hspace{5pt}
        \subfloat[$40^{\mathrm{th}}$ wave period, trough]{\includegraphics[width=.38\linewidth] 
             {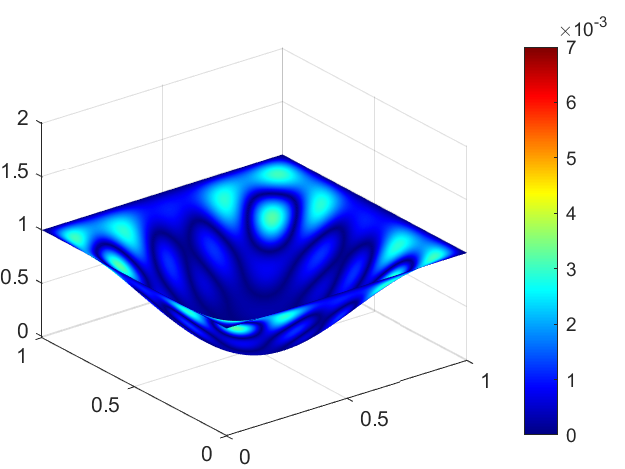}}\\
        \subfloat[$50^{\mathrm{th}}$ wave period, crest]{\includegraphics[width=.38\linewidth] 
             {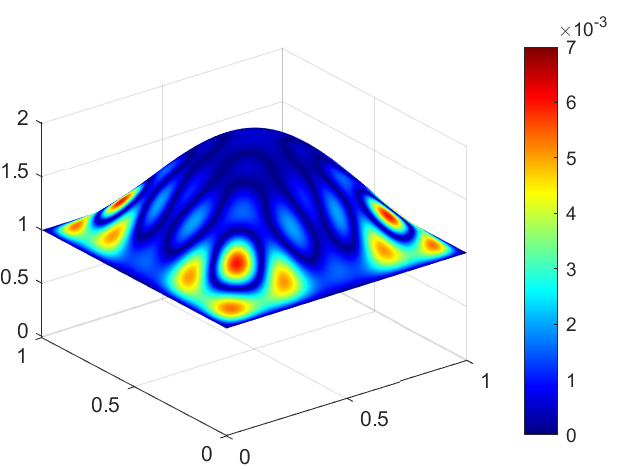}}\hspace{5pt}
        \subfloat[$60^{\mathrm{th}}$ wave period, trough]{\includegraphics[width=.38\linewidth] 
             {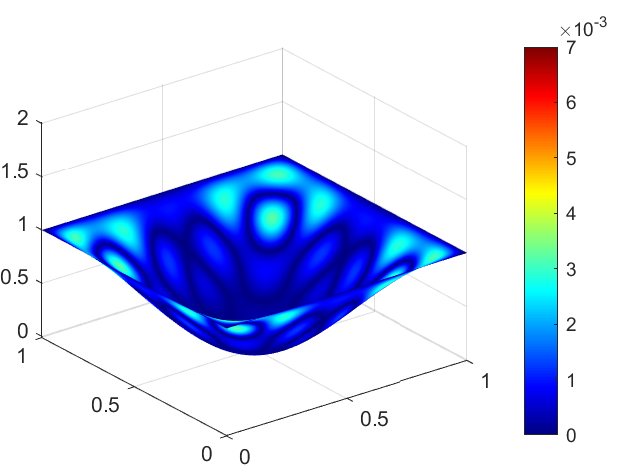}}\\
        \subfloat[$70^{\mathrm{th}}$ wave period, crest]{\includegraphics[width=.38\linewidth] 
             {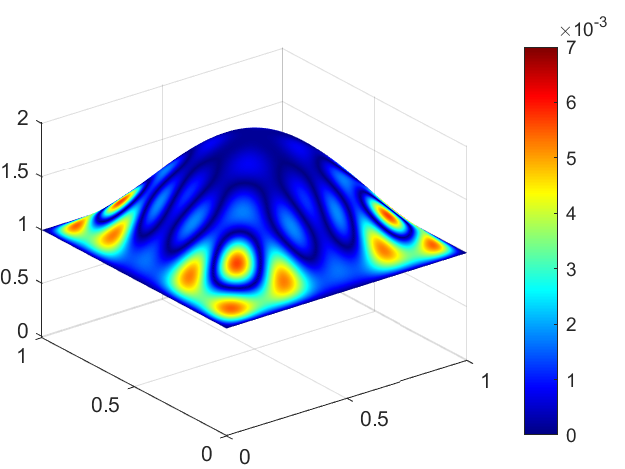}}\hspace{5pt}
        \subfloat[$80^{\mathrm{th}}$ wave period, trough]{\includegraphics[width=.38\linewidth] 
             {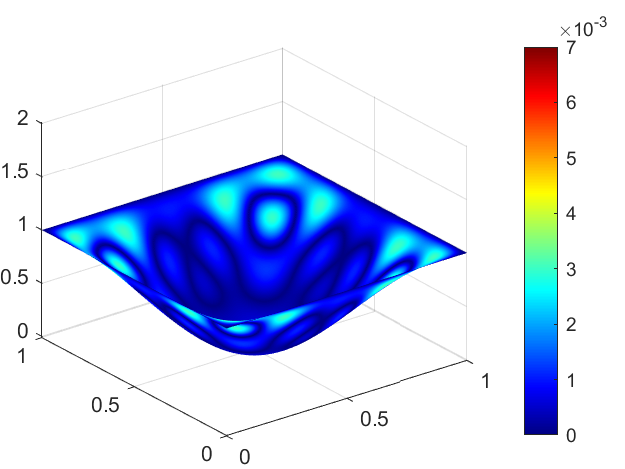}}
    \caption{\textbf{(Example 5 for PDE 1)} Numerical solutions and absolute error distribution of the 2D linear wave equation across various periods using the EC-LS Kansa-CN(AB) method. The left panel displays wave crests at the $10^{\mathrm{th}}$, $30^{\mathrm{th}}$, $50^{\mathrm{th}}$, and $70^{\mathrm{th}}$ wave periods, while the right panel shows wave troughs at the $20^{\mathrm{th}}$, $40^{\mathrm{th}}$, $60^{\mathrm{th}}$, and $80^{\mathrm{th}}$ wave periods.}
    \label{fig:example5 LW2D}
\end{figure}

\begin{figure}[tbhp]
\centering
\begin{overpic}
[width=0.26\textwidth]{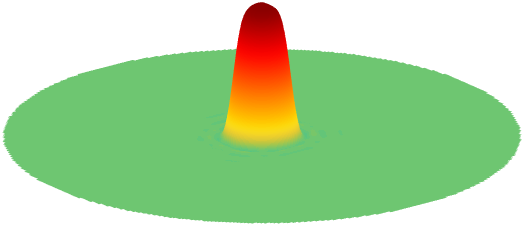}
\put(41,-10){\scriptsize{t = 0 }}
\end{overpic} 
 \hspace{8pt} 
\begin{overpic}
[width=0.26\textwidth]{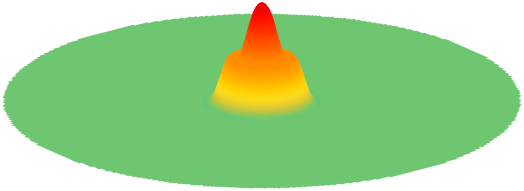}
\put(41,-10){\scriptsize{t = 0.5 }}
\end{overpic} 
\hspace{8pt}
\begin{overpic}
[width=0.26\textwidth]{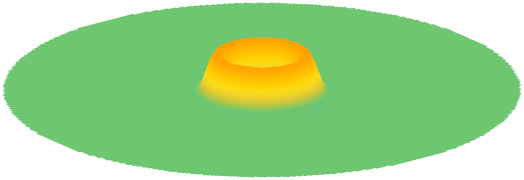}
\put(41,-10){\scriptsize{t = 1 }}
\end{overpic} 
\\
 \vspace{8mm}
\begin{overpic}
[width=0.26\textwidth]{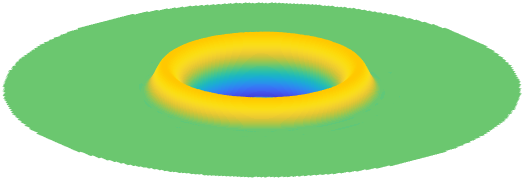}
\put(41,-10){\scriptsize{t = 3 }}
\end{overpic} 
\hspace{8pt}
\begin{overpic}
[width=0.26\textwidth]{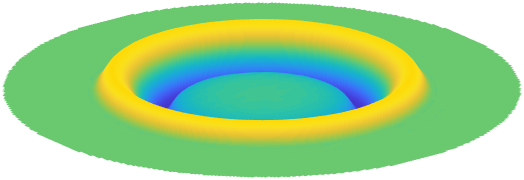}
\put(41,-10){\scriptsize{t = 5 }}
\end{overpic} 
\hspace{5pt}
\begin{overpic}
[width=0.26\textwidth]{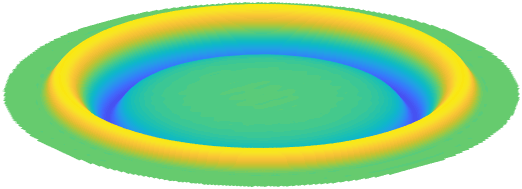}
\put(41,-10){\scriptsize{t = 7 }}
\end{overpic} 
    \caption{\textbf{(Example 5 for PDE 3)} Numerical solutions of EC-LS Kansa-CN method for the nonlinear Klein-Gordon equation with Dirichlet boundary conditions in a circular domain, displaying the temporal evolution of a solitary wave at times \( t = 0.0 ,0.5 , 1.0 , 3.0 , 5.0 , 7.0 \).}
	\label{fig:example5 KG2D Dirichlet}
\end{figure}

\begin{figure}[tbhp]
    \centering
    \subfloat[t = 0,\\ energy = $E_0$ -- 2.56e-13]{\includegraphics[width=.38\linewidth]{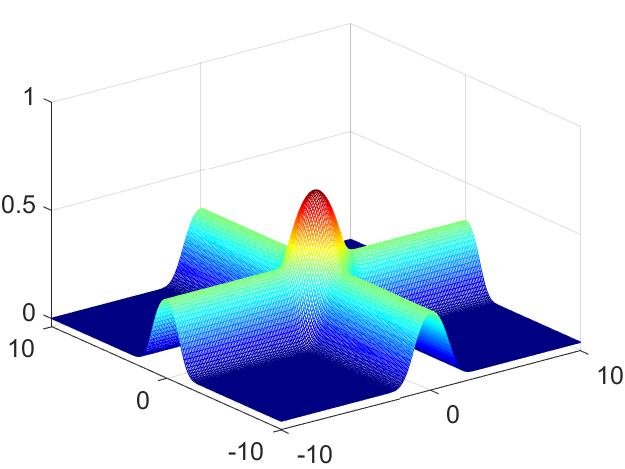}}\hspace{5pt}
    \subfloat[t = 1,\\ energy = $E_0$ + 2.35e-12]{\includegraphics[width=.38\linewidth]{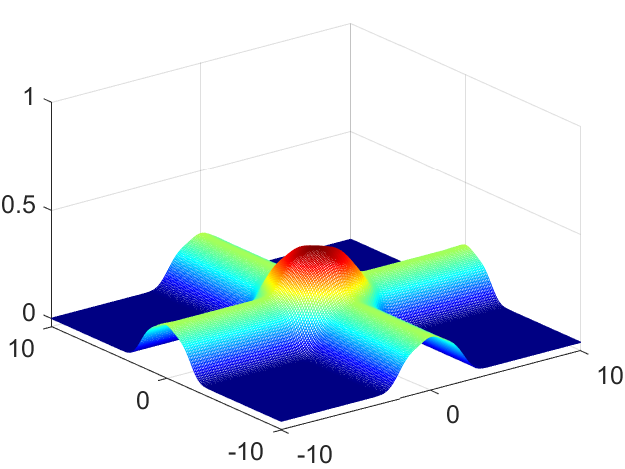}}\\
    \subfloat[t = 4,\\ energy = $E_0$ + 3.67e-12]{\includegraphics[width=.38\linewidth]{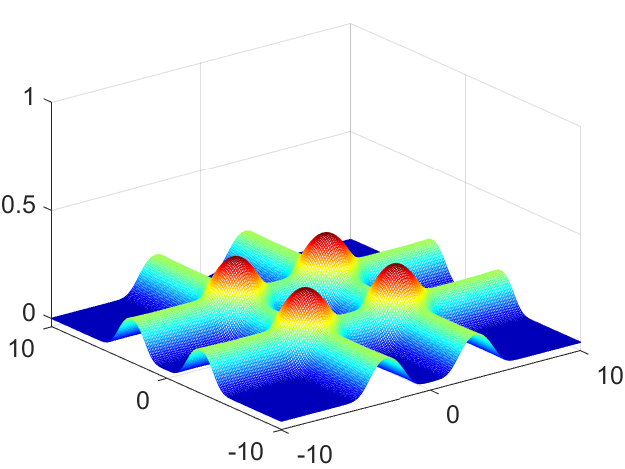}}\hspace{5pt}
    \subfloat[t = 7,\\ energy = $E_0$ -- 4.76e-12]{\includegraphics[width=.38\linewidth]{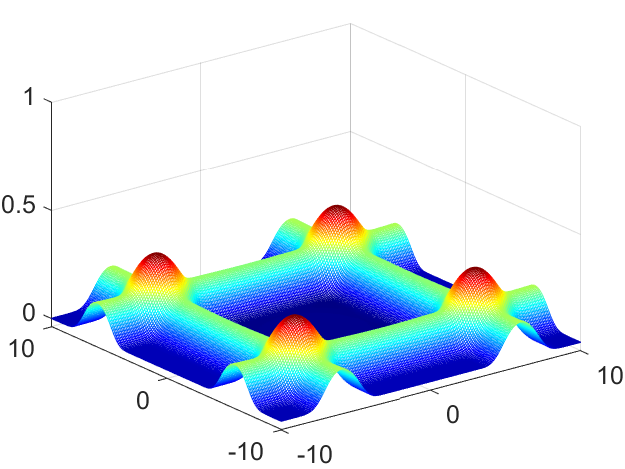}}\\
    \subfloat[t = 10,\\ energy = $E_0$ -- 2.68e-12]{\includegraphics[width=.38\linewidth]{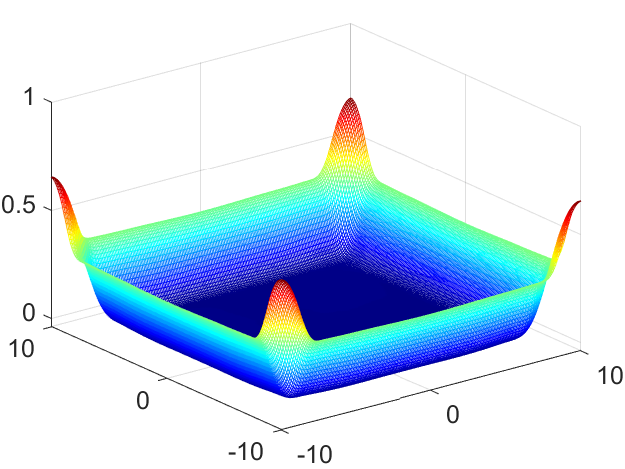}}\hspace{5pt}
    \subfloat[t = 13,\\ energy = $E_0$ + 5.89e-12]{\includegraphics[width=.38\linewidth]{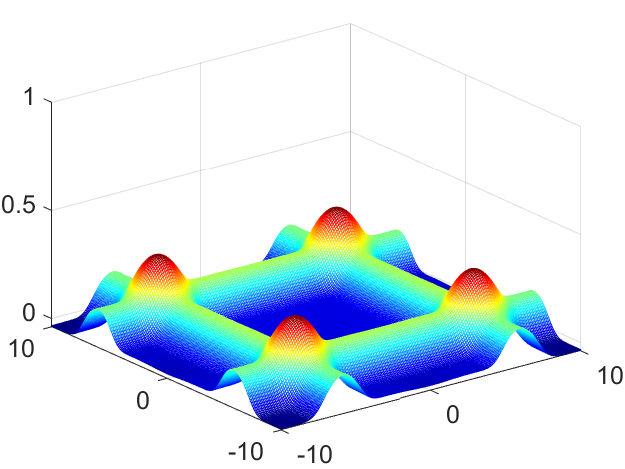}}\\
    \subfloat[t = 17,\\ energy = $E_0$ -- 2.70e-12]{\includegraphics[width=.38\linewidth]{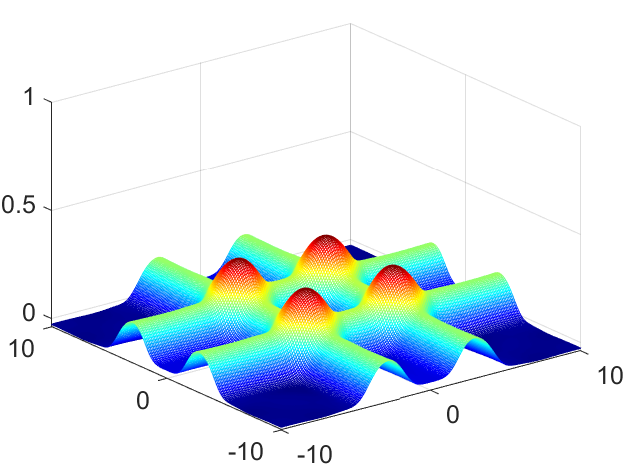}}\hspace{5pt}
    \subfloat[t = 20,\\ energy = $E_0$ -- 1.40e-12]{\includegraphics[width=.38\linewidth]{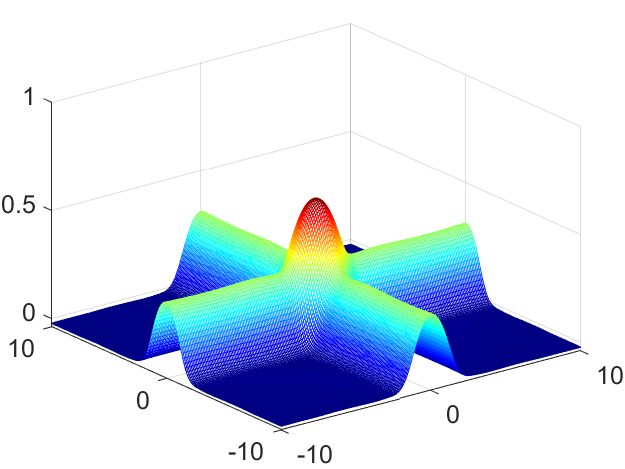}}
\caption{\textbf{(Example 5 for PDE 4)} Time evolution of the 2D nonlinear Klein-Gordon equation $\ddot{u} - \Delta u + u^5 = 0$ solved using EC-LS Kansa-CNAB method. Snapshots are shown for times \( t = 0, 1, 4, 7, 10, 13, 17, 20 \).}
 \label{fig:example5 KG2D Neumann}
\end{figure}

As shown in Figure \ref{fig:example5 LW2D}, the proposed method stably captures the peaks and troughs in the numerical solution of PDE 1, accurately depicting the waveform changes of the periodic wave equation over 80 periods.
Figure \ref{fig:example5 KG2D Dirichlet} describes the propagation process of an isolated wave in a circular area, corresponding to PDE 3.
Figure \ref{fig:example5 KG2D Neumann} reveals the initial soliton splitting into four solitons on opposite sides of the xy-axis, moving outward. When the wave reaches the boundary, it reflects back in the original direction, corresponding to PDE 4.
The methods lacking energy constraint exhibit a significant increase in energy error upon the wave's encounter with the boundary, as illustrated in Figure \ref{fig:example3 different boundary}.
Throughout the simulation, our method consistently demonstrates high energy accuracy. 

\pgfplotsset{compat=1.17}

\section{Conclusion}\label{sec:conclusion}
 In this paper, we introduce a novel meshless energy-conserving method for solving second-order time-dependent Hamiltonian wave equations. 
 The approach combines the discretization of the least-squares Kansa and the time-stepping with a quadratic constraint.
To solve this kind of nonlinear optimization problem, we design a fast iterative algorithm based on GSVD, the Lagrange multiplier method, and
 the Newton method with successive linearization.
 Numerical results demonstrate that our method maintains energy conservation while exhibiting  accuracy and stability. 
 Furthermore, it achieves higher computational precision and efficiency under complex boundary conditions and in long-time simulations. 
\section*{Acknowledgments}
This work was supported by NSFC (No. 12361086, 12001261, 12371379), NSF of Jiangxi Province (No. 20212BAB211020),  NSFC (No. 12101310), NSF of Jiangsu Province (No. BK20210315), the Fundamental Research Funds for the Central Universities (No. 30923010912), the General Research Fund (GRF No. 12301021, 12300922, 12301824) of Hong Kong Research Grant Council, and the NSFC (No. 12201449).
\appendix
\section{Collection of formulas}\label{sec:appendix}
This appendix contains formulas for  implementing the least-squares Kansa methods with the CN and CNAB time-stepping method.
Let $u^0=I_Z \psi_0$ and $v^0=I_Z \psi_1$ be the RBF interpolants of the initial conditions within the working trial space. 

\subsection{Semi-discrete solution in \cref{both eq}}\label{sec:appendix A1}
The formulas for semi-discrete solution $u^k$ in \cref{both eq} are as follows:

\begin{align}
    \tag{Initialization}\label{Init}
    \small{
    u^{1} = \frac{\tau^2}{4} \Delta u^{1} + \left(\frac{\tau^2}{4} \Delta + I\right) u^0 + \tau v^0 - \frac{\tau^2}{2} F'\left(\frac{1}{2}u^{1} + \frac{1}{2} u^0\right), \nonumber
    }
\end{align}

\begin{itemize}
    \item\textbf{CN}: 
    \begin{small}
    \begin{align}
        u^k =& \frac{\tau^2}{4} \Delta u^k + \left(\frac{\tau^2}{2} \Delta + 2I\right) u^{k-1} + \left(\frac{\tau^2}{4} \Delta - I\right) u^{k-2} \nonumber \\
        & - \frac{\tau^2}{2} F'\left(\frac{1}{2}u^k + \frac{1}{2}u^{k-1}\right) - \frac{\tau^2}{2} F'\left(\frac{1}{2}u^{k-1} + \frac{1}{2}u^{k-2}\right),  
         \quad\quad\quad\quad\quad\quad(k \geq 2)\nonumber
    \end{align}
    \end{small}
    \item\textbf{CNAB}: 
    \begin{small}
    \begin{align}
    u^k =& \frac{\tau^2}{4} \Delta u^k + \left(\frac{\tau^2}{4} \Delta + 3I\right) u^{k-1} - 2 u^{k-2} - \tau v^{k-2} - \frac{\tau^2}{2} F'\left(\frac{3}{2} u^{k-1} - \frac{1}{2} u^{k-2}\right),\nonumber
    \\
   &\quad\quad\quad\quad\quad\quad\quad\quad\quad\quad\quad\quad\quad\quad\quad\quad\quad\quad\quad\quad\quad\quad\quad\quad\quad\quad\quad\quad \:\:\:\:(k=2)\nonumber
    \\
    u^k =& \frac{\tau^2}{4} \Delta u^k + \left(\frac{\tau^2}{2} \Delta + 2I\right) u^{k-1} + \left(\frac{\tau^2}{4} \Delta - I\right) u^{k-2} \nonumber \\
        & - \frac{\tau^2}{2} F'\left(\frac{3}{2} u^{k-1} - \frac{1}{2} u^{k-2}\right) - \frac{\tau^2}{2} F'\left(\frac{3}{2} u^{k-2} - \frac{1}{2} u^{k-3}\right).\quad\quad\quad\quad\:( k \geq 3)  \nonumber
    \end{align}
    \end{small}

\end{itemize}

\subsection{Semi-discrete solution in \cref{vk(uk)}}\label{sec:appendix A2}
The formulas for semi-discrete solution $v^k$ in \cref{vk(uk)} are as follows:

\begin{align}
    \tag{Initialization}
    \begin{small}
    v^{1} =\dfrac{2}{\tau}u^1- \dfrac{2}{\tau}u^0-v^0,  \nonumber\quad\quad\quad\quad\quad\quad\quad\quad\quad\quad\quad\quad\quad\quad\quad
    \end{small}
\end{align}

\begin{itemize}
    \item\textbf{CN}: 
    \begin{small}
    \begin{align}
    v^k =& \dfrac{2}{\tau}u^{k}-
    \left(\dfrac{\tau }{4} \Delta+\dfrac{3}{\tau}I\right) u^{k-1}  + \left( - \dfrac{\tau }{4}\Delta+\dfrac{1}{\tau}I\right) u^{k-2}+\dfrac{\tau}{2}
 F^{'}\left(\dfrac{1}{2}u^{k-1}+\dfrac{1}{2}u^{k-2}\right),  \nonumber 
 \\
         &\quad\quad\quad\quad\quad\quad\quad\quad\quad\quad\quad\quad\quad\quad\quad\quad\quad\quad\quad\quad\quad\quad\quad\quad\quad\quad\quad\quad\quad\:\:\,( k \geq 2)  \nonumber
    \end{align}
    \end{small}

    \item\textbf{CNAB}:
    \begin{small}
    \begin{align}
    v^k =& \dfrac{2}{\tau}u^k- \dfrac{4}{\tau}u^{k-1} + \dfrac{2}{\tau}u^{k-2} + v^{k-2},\quad\quad\quad\quad\quad\quad\quad\quad\quad\quad\quad\quad\quad\quad\quad\:\:\:\:(k = 2)
\nonumber\\
    v^k =&\dfrac{2}{\tau}u^{k}-
    \left( \dfrac{\tau }{4} \Delta+\dfrac{3}{\tau}I\right) u^{k-1}  + \left( - \dfrac{\tau }{4}\Delta+\dfrac{1}{\tau}\right) u^{k-2} +\dfrac{\tau}{2}F^{'}\left(\dfrac{3}{2}u^{k-2}-\dfrac{1}{2}u^{k-3}\right).  \nonumber
    \\
     &\quad\quad\quad\quad\quad\quad\quad\quad\quad\quad\quad\quad\quad\quad\quad\quad\quad\quad\quad\quad\quad\quad\quad\quad\quad\quad\quad\quad\quad\:\:\:( k \geq 3)\nonumber
    \end{align}
    \end{small}

\end{itemize}

\subsection{Fully discretized energy in \cref{approx E}}\label{sec:appendix A3}
As for discretizing  \cref{semi energy}, we employ a sufficiently dense set \( P = \{p_1, \dots, p_{n_P}\} \subset \overline{\Omega} \) of quadrature points associated with some quadrature weights \( w_j \).
For any two functions  $f$, $g\in L^2(\Omega)$, their  $L^2$ inner product can be approximated by

\begin{equation}\label{quadrature formula}
\big\langle f, g \big\rangle_{L^2(\Omega)} := \int_\Omega fg \, dx \approx \sum\limits_{j=1}^{n_P} w_j f(p_j) g(p_j) = f(P)^T W g(P),
\end{equation}
where \( W = \text{diag}(w_1, \ldots, w_{n_P}) \) is the \( n_P \times n_P \) diagonal matrix. 

Since \(v^k\) is a linear function of \(u^k\), \((v^k)^2\) can be expressed as a quadratic function of \(u^k\), including interactions with other \(u^j\) for \(j < k\).
Using the quadrature formula \cref{quadrature formula},
the first two terms in the energy in \cref{semi energy}, when fully discretized, are a quadratic function in \( \boldsymbol{\alpha}^k \), whose quadratic form matrix is the  weighted sum of the following two Gram matrices:
\begin{align}
    G_{\Phi}\:\:\: & := 
     \big[\Phi(P, Z)\big]^T W \big[\Phi(P, Z)\big],
     \label{Gram0}\\
    G_{\nabla \Phi} & := 
    \sum_{\substack{\lvert \xi \rvert = 1, \xi \in \mathbb{N}_{0}^{d}}} \big\{ [D^{\xi} \Phi](P, Z) \big\}^T W \big\{ [D^{\xi} \Phi](P, Z) \big\},  
    \label{Gram1}
\end{align}
where \(D^{\xi}\) represents the mixed partial derivatives of function \(\Phi\) according to the multi-index \(\xi\), and
\begin{align*}
[G_{\Phi}]_{ij}\:\:\:&\approx
    {\big\langle \Phi(\cdot, z_i), \Phi(\cdot, z_j) \big\rangle}_{L^2(\Omega)},\\
[G_{\nabla\Phi}]_{ij}&\approx 
     \sum_{\substack{\lvert \xi \rvert = 1, \xi \in \mathbb{N}_{0}^{d}}} {\left\langle D^{\xi} \Phi(\cdot, z_i), D^{\xi} \Phi(\cdot, z_j) \right\rangle}_{L^2(\Omega)}.
\end{align*}
An important fact is the symmetric positive definite (SPD) nature of this quadratic form because Eqs.\eqref{Gram0} and \eqref{Gram1} are SPD.

Along with integrals for lower-order terms and the \( F \)-related final term, we can express the fully discretized energy in terms of RBF coefficients as 
\begin{equation*}
E_w( \boldsymbol{\alpha}^k ) =\mathcal{Q}(\boldsymbol{\alpha}^k) +\mathcal{N}(\boldsymbol{\alpha}^k),
\end{equation*}
where
\begin{align}
\mathcal{Q}(\boldsymbol{\alpha}^k) = &(\boldsymbol\alpha^{k})^{T}\small{\Big(\dfrac{2}{\tau^2} G_{\Phi} + \dfrac{1}{2}G_{\nabla\Phi}\Big)}\boldsymbol\alpha^{k}+
\dfrac{2}{\tau}(\mathbf{d}^{k-1}_{P})^T W [\Phi(P,Z)]\boldsymbol\alpha^{k} + \dfrac{1}{2}
(\mathbf{d}^{k-1}_{P})^T W 
\mathbf{d}^{k-1}_{P}, \label{E_w:Q} \\
\mathcal{N}(\boldsymbol{\alpha}^k) =& \boldsymbol 1_{n_{P}}^{T} W \Big[F\left(\Phi(P,Z) \boldsymbol\alpha^{k}\right)\Big]. \label{E_w:N}
\end{align}

The components $\mathbf{d}^{k-1}_{P}$ of  \( E_w( \boldsymbol{\alpha}^k ) \)  that are independent of the unknown coefficients \(\boldsymbol{\alpha}^k \) are as follows:
\begin{align}
    \tag{Initialization}
    \begin{small}
    \mathbf{d}^{0}_{P} =\Big\{ - \dfrac{2}{\tau}u^0-v^0\Big\}\Big\vert_P,  \qquad\qquad\qquad\qquad\qquad\qquad\qquad\nonumber
    \end{small}
\end{align}

\begin{itemize}
    \item\textbf{CN}:
    \begin{small}
    \begin{align}
    \mathbf{d}^{k-1}_{P} &= -
    \left(\dfrac{\tau }{4} \Delta\Phi(P, Z)+\dfrac{3}{\tau}\Phi(P, Z)\right) \boldsymbol{\alpha}^{k-1}  + \left( - \dfrac{\tau }{4}\Delta\Phi(P, Z)+\dfrac{1}{\tau}\Phi(P, Z)\right) \boldsymbol{\alpha}^{k-2}\nonumber \\
        &\quad+\dfrac{\tau}{2}
 F^{'}\left(\dfrac{1}{2}\Phi(P,Z)\left(\boldsymbol{\alpha}^{k-1} + \boldsymbol{\alpha}^{k-2}\right) \right), 
\nonumber \quad\quad\quad\quad\quad\quad\quad\quad\quad\quad\quad\quad\:(k \geq 2) 
\end{align}
\end{small}
\item\textbf{CNAB}: 
\begin{small}
\begin{align}
    \mathbf{d}^{k-1}_{P} &=- \dfrac{4}{\tau} \Phi(P, Z)\boldsymbol{\alpha}^{k-1} +\dfrac{2}{\tau} \Phi(P, Z)\boldsymbol{\alpha}^{k-2} + v^{k-2}\vert_P, \:\:\quad\quad\quad\quad\quad\quad \quad\quad\:\,\,(k = 2)  
\nonumber
\\
    \mathbf{d}^{k-1}_{P} &=-
    \left( \dfrac{\tau }{4} \Delta\Phi(P, Z)+\dfrac{3}{\tau}\Phi(P, Z)\right) \boldsymbol{\alpha}^{k-1}  + \left( - \dfrac{\tau }{4}\Delta\Phi(P,Z)+\dfrac{1}{\tau}\Phi(P,Z)\right) \boldsymbol{\alpha}^{k-2} \nonumber \\
        &\quad+\dfrac{\tau}{2}
 F^{'}\left(\dfrac{1}{2}\Phi(P,Z)\left(3\boldsymbol{\alpha}^{k-2} - \boldsymbol{\alpha}^{k-3}\right)\right).  \quad\quad\quad\quad\quad\quad\quad\quad\quad\quad\:\:\:\:\quad ( k \geq 3)
\nonumber
\end{align}
\end{small}
\end{itemize}
Additionally, we can analytically compute the initial energy by the initial condition $E_0 := E[\psi_0, \psi_1](0)$ or approximate the initial energy using the RBF interpolants to the initial conditions, i.e., $E_0 \approx E[I_Z\psi_0, I_Z\psi_1](0)$.

\subsection{Matrix forms of optimization problems in Eqs. \eqref{LS sol} and \eqref{full_dis_sol}}\label{sec:appendix A4} 
Components $A$ and $\mathbf{b}_k(\boldsymbol{\eta})$ in the optimization problems Eqs. \eqref{LS sol} and \eqref{full_dis_sol} are:
 \begin{equation}\label{A matrix}
    A=
    \begin{bmatrix}
    -\dfrac{\tau^2}{4}\Delta\Phi(X,Z)+\Phi(X,Z);
    h^{-\theta} [\mathcal{B}\Phi](Y, Z)
    \end{bmatrix},
\end{equation}
and
\begin{equation}\label{b vector}
    \mathbf{b}_k(\boldsymbol{\eta}) =
    \begin{bmatrix}
    \mathbf{b}_{X,k}; h^{-\theta}g(Y)
\end{bmatrix},
\end{equation}
where the \(\mathbf{b}_{X,k}\) in \(\mathbf{b}_k\) are as follows:
\begin{align}
    \tag{Initialization}
    \begin{small}
    \qquad\mathbf{b}_{X,1} =
    \left(\dfrac{\tau^2}{4}\Delta\Phi(X,Z)+\Phi(X,Z)\right)\boldsymbol{\alpha}^0
    +\tau v^0\vert_X-\dfrac{\tau^2}{2} F^{'}\left(\dfrac{1}{2}\Phi(X,Z)(\boldsymbol{\eta} + \boldsymbol{\alpha}^{0})\right),\nonumber
    \end{small}
\end{align}

\begin{itemize}
    \item\textbf{CN}: 
    \begin{small}
    \begin{align}
    \mathbf{b}_{X,k} =&\left(  \dfrac{\tau^2}{2}{\Delta}\Phi(X,Z)+2\Phi(X,Z)  \right)\boldsymbol{\alpha}^{k-1} 
    -\dfrac{\tau^2}{2}F^{'}\left(\dfrac{1}{2}\Phi(X,Z)\left(\boldsymbol\eta + \boldsymbol{\alpha}^{k-1}\right)\right)
     \nonumber \\
        &
        + \left(\dfrac{\tau^2}{4}{\Delta}\Phi(X,Z) -\Phi(X,Z)  \right)\boldsymbol{\alpha}^{k-2}
         -\dfrac{\tau^2}{2} F^{'}\left(\dfrac{1}{2}\Phi(X,Z)\left(\boldsymbol{\alpha}^{k-1} + \boldsymbol{\alpha}^{k-2}\right)\right),
         \nonumber
         \\
         &\quad\quad\quad\quad\quad\quad\quad\quad\quad\quad\quad\quad\quad\quad\quad\quad\quad\quad\quad\quad\quad\quad\quad\quad\quad\quad\quad\quad\quad( k \geq 2)  \nonumber
\end{align}
\end{small}
    \item\textbf{CNAB}: 
    \begin{small}
    \begin{align*}
    \mathbf{b}_{X,k} = &
    \Big(\dfrac{\tau^2}{4}{\Delta}\Phi(X,Z) +3\Phi(X,Z) \Big)\boldsymbol{\alpha}^{k-1}
    -\dfrac{\tau^2}{2}F^{'}\left(\dfrac{1}{2}\Phi(X,Z)\left(3\boldsymbol{\alpha}^{k-1} - \boldsymbol{\alpha}^{k-2}\right)\right)
     \nonumber\\
        &-2\Phi(X,Z) \boldsymbol{\alpha}^{k-2} -\tau v^{k-2}\vert_X,   \quad\quad\quad\quad\quad\quad\quad\quad\quad\quad\quad\quad\quad\quad\quad\quad\:(k = 2)\nonumber
\\
     \mathbf{b}_{X,k} = &
    \Big(\dfrac{\tau^2}{2}{\Delta}\Phi(X,Z) +2\Phi(X,Z) \Big)\boldsymbol{\alpha}^{k-1}
    -\dfrac{\tau^2}{2}F^{'}\left(\dfrac{1}{2}\Phi(X,Z)\left(3\boldsymbol{\alpha}^{k-1} - \boldsymbol{\alpha}^{k-2}\right)\right)
     \nonumber\\
        &+
          \Big(\dfrac{\tau^2}{4}{\Delta}\Phi(X,Z) -\Phi(X,Z) \Big)\boldsymbol{\alpha}^{k-2}
         -\dfrac{\tau^2}{2} F^{'}\left(\dfrac{1}{2}\Phi(X,Z) \left(3\boldsymbol{\alpha}^{k-2} - \boldsymbol{\alpha}^{k-3}\right)\right). \\
          &\quad\quad\quad\quad\quad\quad\quad\quad\quad\quad\quad\quad\quad\quad\quad\quad\quad\quad\quad\quad\quad\quad\quad\quad\quad\quad\quad\quad\quad( k \geq 3) \nonumber
    \end{align*}
    \end{small}
\end{itemize}
The constraint components $B$ and $\mathbf{d}$ from \cref{full_dis_sol} are constructed as:
\begin{equation}\label{B matrix}
    B=
\begin{bmatrix}
\frac{\sqrt{2}}{\tau}W^{\frac{1}{2}} & 0 & \cdots & 0 \\
0 & \frac{\sqrt{2}}{2}  W^{\frac{1}{2}} & \cdots & 0 \\
\vdots & \vdots & \ddots & \vdots \\
0 & 0 & \cdots &\frac{\sqrt{2}}{2}  W^{\frac{1}{2}}
\end{bmatrix}
\begin{bmatrix}
\hphantom{D^{\xi_1}}\Phi(P,Z)\\
D^{\xi_1}\Phi(P,Z)\\
\vdots
\\
D^{\xi_d}\Phi(P,Z)
\end{bmatrix}, 
\end{equation}
and
\begin{equation}\label{d vector}
\mathbf{d}=
\begin{bmatrix}
 -\frac{\sqrt{2}}{2}  W^{\frac{1}{2}}\mathbf{d}^{k-1}_{P}
,
 \boldsymbol 0_{n_P}
,
 \ldots
,
\boldsymbol 0_{n_P}
\end{bmatrix}^T,
\end{equation}
where   
$D^{\xi_i} \Phi  (i = 1, \dots, d)$ denotes the first-order partial derivatives of $\Phi$.
The formula for the nonlinear function in the constraint is as follows:
\begin{equation}\label{N_f}
   \mathcal{N}_{F}(\boldsymbol{\eta}) = \boldsymbol 1_{n_{P}}^{T}W\Big[F\left(\Phi(P,Z)\boldsymbol{\eta}\right)\Big]-E_0. 
\end{equation}



\raggedbottom
\bibliographystyle{elsarticle-num-names}
\bibliography{reference}

\end{document}